\documentclass[reqno]{amsart}
 \newtheorem{thm}{Theorem}[section]

 \theoremstyle{definition}
 
 \theoremstyle{remark}

 \numberwithin{equation}{section}

\usepackage{verbatim}
\usepackage{xcolor}
\usepackage{graphicx}
\usepackage{amssymb, amsmath, amsthm}
\usepackage{amsfonts}
\usepackage{amssymb}
\usepackage{float}
\usepackage{hyperref}
\usepackage{mathrsfs}
\usepackage{enumitem}

\newtheorem{theorem}{Theorem}

\newtheorem{definition}[theorem]{Definition}

\newcommand{\dis}{\displaystyle}

\newcommand{\bK}{\bar{K}_{\infty}}

\newcommand{\dimf}{\text{\rm dim}_{f}}
\newcommand{\A}{{\mathcal A}}
\newcommand{\An}{{\mathcal A}^{\nu}}
\newcommand{\Gg}{\mathcal{G}^{\nu}}
\newcommand{\Ggo}{\mathcal{G}^{\nu_0}}

\newcommand{\K}{K_{\infty}}

\newcommand{\N}{\mathbb N}

\newcommand{\pin}{\pi^{\nu}}
\newcommand{\pino}{\tilde\pi^{\nu}}

\newcommand{\I}{I^*}
\newcommand{\U}{\mathcal U}

\newcommand{\R}{\mathbb R}

\newcommand{\intox}{\int_{\mathbb{T}^d}}

\numberwithin{equation}{section}
\numberwithin{theorem}{section}

\numberwithin{figure}{section}

\begin{document}

%
%
%
%
%
%
%
%
%

\title[Survey on active scalar equations]
 {On a class of forced active scalar equations with small diffusive parameters}

\author{Susan Friedlander}

\address{Department of Mathematics\\
University of Southern California}

\email{susanfri@usc.edu}

\author{Anthony Suen}

\address{Department of Mathematics and Information Technology\\
The Education University of Hong Kong}

\email{acksuen@eduhk.hk}

\dedicatory{In Memory of Louis Nirenberg}

\date{}



\begin{abstract}
Many equations that model fluid behaviour are derived from systems that encompass multiple physical forces. When the equations are written in non dimensional form appropriate to the physics of the situation, the resulting partial differential equations often contain several small parameters. We study a general class of such PDEs called active scalar equations which in specific parameter regimes produce certain well known models for fluid motion. We address various mathematical questions relating to well-posedness, regularity and long time behaviour of the solutions to this general class including vanishing limits of several diffusive parameters.
\end{abstract}

\maketitle
\tableofcontents
\section{Introduction}\label{introduction}

Active scalar equations belong to a class of partial differential equations where the evolution in time of a scalar quantity is governed by the motion of the fluid where the velocity itself varies with this scalar quantity. They have been a topic of considerable study in recent years, in part because they arise in many physical models and in part because they present challenging nonlinear PDEs. Various active scalar equations, such as surface quasi-geostrophic equation (SQG) and drift-diffusion equations, have received considerable attention in the past decade because of the challenging nature of the delicate balance between the nonlinear term and the dissipative term. The physics of an active scalar equation is encoded in the constitutive law that relates the transport velocity vector $u$ with a scalar field $\theta$. This law produces a differential operator that when applied to the scalar field determines the velocity. The singular or smoothing properties of this operator are closely connected with the mathematics of the nonlinear advection equation for $\theta$.

More precisely, we are interested in an abstract class of active scalar equations in $\mathbb{T}^d\times(0,\infty)=[0,2\pi]^d\times(0,\infty)$ with $d\in\{2,3\}$ \footnote{We point out that, most of the results given in our work hold for $d\ge2$.} of the following form
\begin{align}
\label{abstract active scalar eqn} \left\{ \begin{array}{l}
\partial_t\theta+u\cdot\nabla\theta=-\kappa\Lambda^{\gamma}\theta+S, \\
u_j[\theta]=\partial_{x_i} T_{ij}^{\nu}[\theta],\theta(x,0)=\theta_0(x)
\end{array}\right.
\end{align}
where $\nu\ge0$, $\kappa\ge0$, $\gamma\in(0,2]$ and $\Lambda:=\sqrt{-\Delta}$. Here $\theta_0$ is the initial datum and $S=S(x)$ is a given function that represents the forcing of the system. We assume that \footnote{Such mean zero assumption is common in many physical models which include SQG equation and MG equation; see \cite{CTV14a} and \cite{FS18} for example.}
\begin{equation}\label{zero mean assumption on data and forcing}
\int_{\mathbb{T}^d}\theta_0(x)dx=\int_{\mathbb{T}^d}S(x)=0
\end{equation}
which immediately implies that $\theta$ obeys
\begin{equation}\label{zero mean assumption}
\int_{\mathbb{T}^d}\theta(x,t)dx=0,\qquad\forall t\ge0. 
\end{equation}
$\{T_{ij}^{\nu}\}_{\nu\ge0}$ is a sequence of operators which satisfy:
\begin{enumerate}
\item[(A1)] For all $\nu\ge0$, $\partial_i\partial_j T^{\nu}_{ij}f=0$ for any smooth functions $f$.
\item[(A2)] There exists a constant $C>0$ independent of $\nu$, such that for all $i,j\in\{1,\dots,d\}$, $$\sup_{\nu\in(0,1]}\sup_{\{k\in\mathbb{Z}^3\}}|\widehat T^{\nu}_{ij}(k)|\le C\mbox{; }\sup_{\{k\in\mathbb{Z}^3\}}|\widehat T^0_{ij}(k)|\le C, \mbox{ where } T^0_{ij}=T_{ij}^{\nu}\Big|_{\nu=0}.$$
\item[(A3)] For each $\nu>0$, there exists a constant $C_{\nu}>0$ such that for all $1\le i,j\le d$, $$|\widehat T^{\nu}_{ij}(k)|\le C_{\nu}|k|^{-3}, \forall k\in\mathbb{Z}^d.$$
\item[(A4)] $T_{ij}^\nu:L^\infty\rightarrow BMO$ are bounded operators for all $\nu\ge0$.
\item[(A5)] For each $1\le j\le d$ and $g\in L^2$, \qquad $\dis\lim_{\nu\rightarrow0}\sum_{k\in\mathbb{Z}^3}|\widehat{T^{\nu}_{ij}}(k)-\widehat{T^0_{ij}}(k)|^2|\widehat{g}(k)|^2=0.$
\end{enumerate}

Our motivation for addressing such a class of active scalar equations mainly comes from several different physical systems, all of them take the form \eqref{abstract active scalar eqn} under particular parameter regimes:
\begin{itemize}
\item[1.] The first example comes from a model proposed by Moffatt and Loper \cite{ML94}, Moffatt \cite{M78} for magenetostrophic turbulence in the Earth's fluid core. Under the postulates in \cite{ML94}, the governing equation becomes a 3D active scalar equation for a temperature field $\theta$
\begin{align}
\label{MG equation} \left\{ \begin{array}{l}
\partial_t \theta + u \cdot \nabla \theta = \kappa \Delta \theta + S,\\
u=M^{\nu}[\theta],\theta(x,0)=\theta_0(x).
\end{array}\right.
\end{align}
The expressions for the Fourier multiplier symbol $\widehat{M}^\nu$ are explicitly given by
\begin{align} 
\label{MG Fourier symbol} \left\{ \begin{array}{l}
\widehat M^{\nu}_1(k)=[k_2k_3|k|^2-k_1k_3(k_2^2+\nu|k|^4)]D(k)^{-1},\\
\widehat M^{\nu}_2(k)=[-k_1k_3|k|^2-k_2k_3(k_2^2+\nu|k|^4)]D(k)^{-1},\\
\widehat M^{\nu}_3(k)=[(k_1^2+k_2^2)(k_2^2+\nu|k|^4)]D(k)^{-1},
\end{array}\right.
\end{align}
where $k = (k_1, k_2, k_3) \in \mathbb{Z}^3$ is the Fourier variable and $D(k)=|k|^2k_3^2+(k_2^2+\nu|k|^4)^2$. The nonlinear equation \eqref{MG equation} with $u$ related to $\theta$ via the operator $M^{\nu}$ is known as the magnetogeostrophic (MG$^{\nu}$) equation (or simply MG equation) and its mathematical properties have been addressed in a series of papers which include  \cite{FS15}, \cite{FS18}, \cite{FS19}, \cite{FS20}, \cite{FV11a}, \cite{FV11b}, \cite{FRV12}, \cite{FRV14}. The behaviour of the MG$^{\nu}$ equation is strikingly different when the parameters $\nu$ and $\kappa$ are present (i.e. positive) or absent (i.e. zero). As the Fourier multiplier symbols $\widehat{M}^0$ given by \eqref{MG Fourier symbol} with $\nu = 0$ are {\it not} bounded in all regions of Fourier space \cite{FV11a}, when $\nu = 0$ the relation between $u$ and $\theta$ is given by a {\it singular} operator of order 1. The implications of such fact for the inviscid MG$^{0}$ equation are summarized in the survey article by Friedlander, Rusin and Vicol \cite{FRV12}. In particular, when $\kappa > 0$ the inviscid but thermally dissipative MG$^{0}$ equation is globally well-possed \cite{FV11a}; on the other hand when $\nu = 0$ {\it and} $\kappa = 0$, the singular inviscid MG$^{0}$ equation is {\it ill-possed} in the sense of Hadamard in any Sobolev space \cite{FV11b}. In \cite{FS15}, Friedlander and Suen first addressed the system \eqref{MG equation}-\eqref{MG Fourier symbol} for $\nu>0$ and obtained well-posedness results in Sobolev space. In a series of papers \cite{FS18}-\cite{FS19} the authors further examined the limit of vanishing viscosity as $\nu\to0$ in the case when $\kappa > 0$ and $\kappa = 0$ and the long-time behaviour of solutions. They proved global existence of classical solutions to the forced MG$^{\nu}$ equations and obtained convergences of solutions as $\nu$ vanishes. Moreover, it was shown in \cite{FS20} that the equations \eqref{MG equation}-\eqref{MG Fourier symbol} possess global attractors with various interesting properties.
\item[2.] The second example of a physical system comes from incompressible flow in a porous medium. It can be modeled by an active scalar equation where a small smoothing parameter enters into the constitutive law. Different from the usual incompressible porous media (IPM) equation, the incompressible porous media Brinkman (IPMB$^\nu$) equation with an ``effective viscosity'' $\nu$ is derived via a modified Darcy's Law as suggested by Brinkman \cite{B49}. The IPM equation becomes the limiting case for IPMB$^\nu$ equation when $\nu=0$. The 2D equation relating the velocity $u$, the density $\theta$ and the pressure $P$ is given in non-dimensional form by
\begin{align}
u &= - \nabla P - e_2 \theta + \nu \Delta u \label{IPMB equation linear system 1}\\
\nabla \cdot u &= 0 \label{IPMB equation linear system 2}
\end{align}
which produces the constitutive law
\begin{align}\label{IPMB constitutive law}
u &= (1 - \nu \Delta)^{-1} [-\nabla \cdot (-\Delta)^{-1} e_2 \cdot \nabla \theta - e_2\theta]\notag\\
&= (1-\nu \Delta)^{-1} R^{\bot} R_1 \theta=M^\nu[\theta]
\end{align}
where $R = (R_1, R_2)$ is the vector of Riesz transforms and $e_2=(0,1)$. The corresponding active scalar equation is thus given by
\begin{align}
\label{active scalar equation general IPMB}
\left\{ \begin{array}{l}
\partial_t\theta^\nu+(u^\nu\cdot\nabla)\theta^\nu=0, \\
u^\nu=M^{\nu}[\theta^{\nu}],\theta^\nu(x,0)=\theta_0(x),
\end{array}\right.
\end{align}
where the 2D components of the Fourier multiplier symbol corresponding to $\widehat{M}^\nu$ as in \eqref{active scalar equation general IPMB} are
\begin{equation}\label{Fourier multiplier symbol IPMB}
\frac{1}{1 + \nu (k^2_1 + k^2_2)} \left(\frac{k_1k_2}{k^2_1 + k^2_2} \quad , \quad \frac{-k^2_1}{(k^2_1 + k^2_2)}\right)
\end{equation}
Similar to the case of MG$^{\kappa,\nu}$ equation, there is a noticeable difference in the operator $M^\nu$ between the two cases $\nu > 0$ and $\nu = 0$: the operator is smoothing of order 2 when $\nu>0$, while for $\nu=0$ the operator is singular of order zero.

The well known IPM equations, i.e. \eqref{IPMB constitutive law}-\eqref{Fourier multiplier symbol IPMB} without the effective viscosity $\nu$, have been studied in a number of papers (see \cite{CFG11}, \cite{CGO07} and the reference therein). As we pointed out before, when $\nu = 0$ the operator in \eqref{IPMB constitutive law} becomes a singular integral operator of order zero, which is similar to the case for the SQG equations. Yet there is a crucial difference between the two operators: the SQG operator is {\it odd} while the IPM operator is {\it even}. Implications for well/ill posedness due to the odd/even structure of the operator in an active scalar equation are further explored in \cite{FGWV12}, \cite{FRV14}, \cite{KVW16}. In a recent work, Friedlander and Suen \cite{FS19} studied the system \eqref{IPMB constitutive law}-\eqref{Fourier multiplier symbol IPMB} in the limit of vanishing viscosity as $\nu\to0$ and obtained convergence results in Sobolev spaces, which are {\it not} available for the case of MG equation when $\kappa=0$. The foremost difference is that the MG$^{0}$ operator is singular of order 1 where as the IPM operator is singular of order zero. In view of such difference, for the “smoother" IPM case the convergence results are valid in Sobolev spaces rather than the Gevrey-class convergence results for the MG equation.
\item[3.] The third physical example comes from the {\it modified} surface quasi-geostrophic (SQG$^{\kappa,\nu}$) equation. It relates the potential temperature $\theta^{\kappa,\nu}$ and the flow velocity $u^{\kappa,\nu}$ as follows:
\begin{align}
\label{surface quasi-geostrophic general SQG}
\left\{ \begin{array}{l}
\partial_t\theta^{\kappa,\nu}+u^{\kappa,\nu}\cdot\nabla\theta^{\kappa,\nu}=-\kappa(-\Delta)^\gamma\theta^{\kappa,\nu}+S, \\
u^{\kappa,\nu}=M^{\nu}[\theta^{\kappa,\nu}],\theta^{\kappa,\nu}(x,0)=\theta_0(x),
\end{array}\right.
\end{align} 
where the Fourier multiplier symbol for $\widehat{M}^\nu$ is given by
\begin{align}\label{Fourier multiplier symbol SQG}
\frac{1}{1 + \nu (k^2_1 + k^2_2)} \left(\frac{k_2}{\sqrt{k^2_1 + k^2_2}} \quad , \quad \frac{-k_1}{\sqrt{k^2_1 + k^2_2}}\right).
\end{align}
When $\nu=0$, the system \eqref{surface quasi-geostrophic general SQG}-\eqref{Fourier multiplier symbol SQG} reduces to the well-known SQG equation which has been investigated by many researchers \cite{BSV16}, \cite{CC10}, \cite{CCCGW12}, \cite{CCV16}, \cite{CI17}, \cite{CTV14b}, \cite{CV16}, \cite{W0405}. For $\kappa>0$ and $\gamma=\frac{1}{2}$, it can be used as a model for studying the temperature distribution $\theta$ on the 2D boundary of a rapidly rotating fluid with small Rossby and Ekman numbers \cite{CW99}, while for $\kappa=0$, the inviscid model can be  applied for studying frontogenesis in meteorology \cite{La17}. It is also worth mentioning that the 3D analog of \eqref{surface quasi-geostrophic general SQG}-\eqref{Fourier multiplier symbol SQG} is widely used as a testing model for the vorticity evolution of the 3D Navier-Stokes equations \cite{CCCGW12}.   
\end{itemize}

We point out that, in view of the physical models as mentioned above, conditions (A1)--(A5) are both mathematically and physically important for studying the active scalar equations \eqref{abstract active scalar eqn}-\eqref{zero mean assumption}, which can be explained as follows:

\begin{itemize}
\item Condition (A1) implies the drift velocity $u$ in \eqref{abstract active scalar eqn} is {\it divergence-free}, which is compatible with the {\it incompressibility} of the fluid described by those physical models. On the other hand, condition (A2) requires a {\it uniform} bound (independent of $\nu$) on the Fourier multiplier symbol $\widehat T^{\nu}$. This conditions implies that the operator $\partial_{x_i} T_{ij}^{\nu}$ is at most {\it singular} of order 1 for $\nu\ge0$, which is consistent with the cases for $\widehat M^\nu$ given in the MG$^{\kappa,\nu}$ equation. We remark here that, however, for the cases of IPMB$^{\nu}$ and SQG$^{\kappa,\nu}$ equation, one can replace the condition (A2) by the following condition:
\begin{itemize}
\item[(A2$^*$)] There exists a constant $C_0>0$ independent of $\nu$, such that for all $i,j\in\{1,\dots,d\}$, $$\sup_{\nu\in(0,1]}\sup_{\{k\in\mathbb{Z}^3\}}|\widehat {\partial_{x_i}T^{\nu}_{ij}}(k)|\le C_0\mbox{; }\sup_{\{k\in\mathbb{Z}^3\}}|\widehat {\partial_{x_i}T^0_{ij}}(k)|\le C_0, \mbox{ where } T^0_{ij}=T_{ij}^{\nu}\Big|_{\nu=0}.$$
\end{itemize}
Condition (A2$^*$) requires that the operators $\partial_{x_i}T^{\nu}_{ij}$ are {\it less} singular than those given by condition (A2), which allows us to obtain better regularity results on IPMB$^{\nu}$ and SQG$^{\kappa,\nu}$ equation. 
\item Condition (A4) imposes a minimal regularity requirement on the operators $T_{ij}^\nu$, while condition (A3) describes the {\it smoothing effect} given by the parameter $\nu$. Condition (A3) implies that the operators $\partial_{x_i} T_{ij}^{\nu}$ are smoothing of degree two for $\nu>0$, which is crucial for proving well-posedness for the system \eqref{abstract active scalar eqn} especially when $\kappa=0$ (refer to \cite{FS18} and \cite{FV11a} for the striking differences between the cases $\nu>0$ and $\nu=0$). 
\item For condition (A5), we notice that the Fourier multiplier symbols $\widehat M^\nu$ given in either \eqref{MG Fourier symbol}, \eqref{Fourier multiplier symbol IPMB} or \eqref{Fourier multiplier symbol SQG} do satisfy condition (A5), which implies that the drift velocity $u$ converges to $u^{\kappa,0}$ strongly in $L^2$. Such behaviour of $T_{ij}^{\nu}$ allows strong convergence of solutions of \eqref{abstract active scalar eqn} (refer to Theorem~7.4 in \cite{FS18}), despite the fact that $T_{ij}^{\nu}$ are all matrices of zero-order pseudo-differential operators as required by condition (A4).
\end{itemize}

In section~\ref{s:Diffusive active scalar equations}, we state the results proved in \cite{FS18} for the diffusive active scalar equations, which are the equations \eqref{abstract active scalar eqn}-\eqref{zero mean assumption} for $\kappa>0$, $\nu\ge0$ and $\gamma=2$. In this case, the system \eqref{abstract active scalar eqn}-\eqref{zero mean assumption} is globally well-posed. To prove this fact, we control the term $\|\theta(t)\|_{L^\infty}$ which can be done by using De Giorgi techniques. Having established the global well-posedness of \eqref{abstract active scalar eqn}-\eqref{zero mean assumption}, we proceed to address the convergence of solutions as $\nu\to0$ which is based on some uniform estimates on $\theta$ which are independent of $\nu$. Moreover, we define a weak solution to the MG$^0$ equation which we call a ``{\it vanishing viscosity}''{\it solution} and prove the existence of a compact global attractor in $L^2(\mathbb{T}^3)$ for the MG$^\nu$ equations \eqref{MG equation}-\eqref{MG Fourier symbol}. We further obtain the upper semicontinuity of the global attractor as $\nu$ vanishes. 

In section~\ref{s:Non-diffusive active scalar equations}, we examine the non-diffusive active scalar equations \eqref{abstract active scalar eqn}-\eqref{zero mean assumption} by considering $\kappa=0$ and $\nu\ge0$ (cf \cite{FS19}). For the case of $\nu>0$, the operators $\partial_x T^\nu$ are smoothing order 2 which give rise to well-posedness for \eqref{abstract active scalar eqn}-\eqref{zero mean assumption} with $\kappa=0$, $\nu>0$ and $\theta_0\in W^{s,d}$ for $s\ge0$ and smooth forcing term $S$, and the results will be discussed in subsection~\ref{s:Well-posedness in various spaces nu+ non-diffusive}. On the other hand, for the case when $\kappa=0$ and $\nu=0$, in general the system \eqref{abstract active scalar eqn}-\eqref{zero mean assumption} fails to be well-posed in Sobolev spaces. In \cite{FV11b}, It was proved that the singular inviscid MG$^0$ equation is {\it ill-possed} in the sense of Hadamard in any Sobolev space. Yet it is possible to obtain the local existence and uniqueness of solutions to \eqref{abstract active scalar eqn}-\eqref{zero mean assumption} with $\kappa=\nu=0$ in spaces of real-analytic functions, owing to the fact that the derivative loss in the nonlinearity $u\cdot\nabla\theta$ is of order at most one. We thus prove the local-in-time existence of Gevrey-class solutions which are summarised in subsection~\ref{s:Well-posedness in various spaces nu0 non-diffusive}. We further address the convergence of Gevrey-class solutions as $\nu\to0$ and apply the claimed results to inviscid MG$^\nu$ equation \eqref{MG equation}-\eqref{MG Fourier symbol} with $\kappa=0$ and IPMB$^\nu$ equation \eqref{IPMB constitutive law}-\eqref{Fourier multiplier symbol IPMB}.

In section~\ref{s:Fractionally diffusive active scalar equations}, we discuss the result obtained in \cite{FS20} by investigating the properties of \eqref{abstract active scalar eqn}-\eqref{zero mean assumption} in the full range $\gamma\in[0,2]$. More precisely we prove the existence and convergence of solutions in various space for the cases $\nu>0$ and $\nu=0$, which are applicable for the critical SQG equation \eqref{surface quasi-geostrophic general SQG}-\eqref{Fourier multiplier symbol SQG} with $\nu=0$. We then address the long-time behaviour for solutions when $\kappa$, $\nu>0$ and obtain global attractors in $H^1$. We further prove some additional properties of the attractors which will be given in subsection~\ref{s:Long time behaviour for solutions when nu > 0 and kappa > 0}. The results on the global attractors can be applied to MG$^\nu$ equation \eqref{MG equation}-\eqref{MG Fourier symbol} which are related to those discussed in section~\ref{s:Diffusive active scalar equations}.

\section{Diffusive active scalar equations}\label{s:Diffusive active scalar equations}

In this section, we discuss the global existence of classical solutions to \eqref{abstract active scalar eqn}-\eqref{zero mean assumption} for $\kappa>0$, $\nu\ge0$ and $\gamma=2$, and study the convergence of solutions as $\nu$ vanishes. Specifically, we consider the following abstract system:
\begin{align}
\label{abstract active scalar eqn diffusive} \left\{ \begin{array}{l}
\partial_t\theta+u\cdot\nabla\theta=\kappa\Delta\theta+S, \\
u_j[\theta]=\partial_{x_i} T_{ij}^{\nu}[\theta],\theta(x,0)=\theta_0(x)
\end{array}\right.
\end{align}
with $\nu\ge0$, and $S$, $\theta_0$ satisfy the zero mean assumptions \eqref{zero mean assumption on data and forcing}-\eqref{zero mean assumption}. The results can then be applied to MG$^\nu$ equations, which allows us to obtain the existence of a compact global attractor $\{\mathcal{A}^\nu\}_{\nu\ge0}$ in $L^2(\mathbb{T}^3)$ including the critical equation where $\nu=0$. 

\subsection{Existence of smooth solutions and convergence as $\nu\to0$}\label{s:Existence of smooth solutions and convergence diffusive}

Friedlander and Vicol \cite{FV11a} analyzed the unforced $S=0$ system \eqref{abstract active scalar eqn diffusive} with the viscosity parameter $\nu$ set to zero, i.e. the unforced MG$^0$ equation. In this situation the drift diffusion equation \eqref{abstract active scalar eqn diffusive} is critical in the sense of the derivative balance between the advection and the diffusion term. They used De Giorgi techniques to obtain global well-posedness results for the unforced critical MG$^0$ equation in a similar manner to the proof of global well-possedness given by Caffarelli and Vaseur \cite{CV10} for the critical SQG equation. Following this work, we verify that the technical details of the De Giorgi techniques are, in fact, valid for drift diffusion equations with a smooth force. More precisely, the result is given by

\begin{thm}[Friedlander and Suen \cite{FS18}]\label{Existence of smooth solutions full diffusive thm}
Let $\theta_0\in L^2$, $S\in C^{\infty}$ and $\kappa>0$ be given, and assume that $\{T_{ij}^{\nu}\}_{\nu\ge0}$ satisfy conditions {\rm (A1)--(A5)}. For $\gamma=2$, there exists a classical solution $\theta^{\nu}(t,x)\in C^{\infty}((0,\infty)\times \mathbb{T}^d)$ of \eqref{abstract active scalar eqn diffusive}, evolving from $\theta_0$ for all $\nu\ge0$.
\end{thm}

Furthermore, for any $\nu\ge0$, we prove that the smooth solutions $\theta^{\nu}$ obtained in Theorem~\ref{Existence of smooth solutions full diffusive thm} satisfy a uniform bound which is independent of $\nu$:

\begin{thm}[Friedlander and Suen \cite{FS18}]\label{uniform bound on solution diffusive thm}
Assume that the hypotheses and notations of Theorem~\ref{Existence of smooth solutions full diffusive thm} are in force. Then given $0<t_1<t_2$ and $s\ge0$, there exists a positive constant $C$ which depends on $C_{0}$, $t_1$, $t_2$, $s$, $d$, $\kappa$, $S$, $\|\theta_0\|_{L^2}$ but independent of $\nu$ such that
\begin{align}\label{uniform bound on solution diffusive}
\sup_{t\in[t_1,t_2]}\|\theta^{\nu}(t,\cdot)\|_{ H^s}+\int_{t_1}^{t_2}\|\theta^{\nu}(t,\cdot)\|^2_{H^{s+1}}dt\le C(C_{0},t_1,t_2,s,d,\kappa,S,\|\theta_0\|_{L^2}),
\end{align}
where $C_{0}>0$ is the constant as stated in condition {\rm (A5)}.
\end{thm}

Using Theorem~\ref{uniform bound on solution diffusive thm}, we obtain the convergence of $\theta^{\nu}$ as $\nu\rightarrow0$:

\begin{thm}[Friedlander and Suen \cite{FS18}]\label{Convergence of solutions diffusive thm}
Assume that the hypotheses and notations of Theorem~\ref{Existence of smooth solutions full diffusive thm} are in force. For $\gamma=2$, if $\theta^{\nu},\theta$ are $C^\infty$ smooth classical solutions of the system \eqref{abstract active scalar eqn diffusive} for $\nu>0$ and $\nu=0$ respectively with initial data $\theta_0$, then given $\tau>0$, for all $s\ge0$, we have
\setcounter{equation}{14}
\begin{align}\label{1.6-}
\lim_{\nu\rightarrow0} \|(\theta^{\nu}-\theta)(t,\cdot)\|_{ H^s}=0,
\end{align}
whenever $t\ge\tau$.
\end{thm}

\subsection{The MG equations and existence of compact global attractor}\label{s:The MG equations and existence of compact global attractor diffusive}

We apply the results obtained in subsection~\ref{s:Existence of smooth solutions and convergence diffusive} to the case for MG equations given by \eqref{MG equation}-\eqref{MG Fourier symbol}. We write 
\begin{align*}
u^{\nu}_j =M^{\nu}_j[\theta^{\nu}]:=\partial_{i} T_{ij}^{\nu},
\end{align*}
where we have denoted 
\begin{align*}
\mbox{$T_{ij}^{\nu}:=-\partial_{i}(-\Delta)^{-1}M^{\nu}_j$ for $\nu\ge0$}
\end{align*}
and $M^{\nu}$ is defined by the inverse Fourier transform of \eqref{MG Fourier symbol}. Using \cite[Lemma~5.1]{FS18}, one can verify that there are constants $C_1,C_2>0$ independent of $\nu$ such that, for all $1\le i,j\le 3$, 
\begin{align*}
\sup_{\nu\in(0,1]}\sup_{\{k\in\mathbb{Z}^3:k\neq0\}}|\widehat T^{\nu}_{ij}(k)|\le\sup_{\nu\in(0,1]}\sup_{\{k\in\mathbb{Z}^3:k\neq0\}}\frac{|\widehat M^{\nu}(k)|}{|k|}\le C_1,
\end{align*}
\begin{align*}
\sup_{\{k\in\mathbb{Z}^3:k\neq0\}}|\widehat T^0_{ij}(k)|\le\sup_{\{k\in\mathbb{Z}^3:k\neq0\}}\frac{|\widehat M^0(k)|}{|k|}\le C_2.
\end{align*}
Hence conditions (A1)--(A5) are satisfied. Theorem~\ref{Existence of smooth solutions full diffusive thm} and Theorem~\ref{uniform bound on solution diffusive thm} can therefore be applied to the MG equations in order to obtain the global-in-time existence and convergence of smooth solutions:

\begin{thm}[Friedlander and Suen \cite{FS18}]\label{existence MG thm}
Let $\theta_0\in L^2$, $S\in C^{\infty}$ and $\kappa>0$ be given. There exists a classical solution $\theta^{\nu}(t,x)\in C^{\infty}((0,\infty)\times \mathbb{T}^3)$ of \eqref{MG equation}-\eqref{MG Fourier symbol}, evolving from $\theta_0$ for all $\nu\ge0$.
\end{thm}

\begin{thm}[Friedlander and Suen \cite{FS18}]\label{convergence MG thm}
Let $\theta_0\in L^2$, $S\in C^{\infty}$ and $\kappa>0$ be given. Then if $\theta^{\nu},\theta$ are $C^\infty$ smooth classical solutions of \eqref{MG equation}-\eqref{MG Fourier symbol} for $\nu>0$ and $\nu=0$ respectively with initial data $\theta_0$, then given $\tau>0$, for all $s\ge0$, we have
\begin{align*}
\lim_{\nu\rightarrow0} \|(\theta^{\nu}-\theta)(t,\cdot)\|_{ H^s}=0,
\end{align*} 
whenever $t\ge\tau$.
\end{thm}

With the results of Theorems~\ref{existence MG thm} and Theorem~\ref{convergence MG thm} in place, we define a weak solution to the MG$^0$ equation which we call a ``{\it vanishing viscosity}''{\it solution}. We use this concept to prove the existence of a compact global attractor in $L^2(\mathbb{T}^3)$ for the MG$^\nu$ equations \eqref{MG equation}-\eqref{MG Fourier symbol} including the critical equation where $\nu=0$, and we further obtain the upper semicontinuity of the global attractor as $\nu$ vanishes. 

\begin{definition}
A {\it weak solution} to \eqref{MG equation}-\eqref{MG Fourier symbol} with $\nu=0$ is a function $\theta\in C_w([0,T];L^2(\mathbb{T}^3))$ with zero spatial mean that satisfies \eqref{MG equation} in a distributional sense. That is, for any $\phi\in C_0^\infty((0,T)\times\mathbb{T}^3)$,
\begin{align*}
-\int_0^T\langle\theta,\phi_t\rangle dt-\int_0^T\langle u\theta, \nabla\phi\rangle dt+\kappa\int_0^T\langle\nabla\theta,\nabla\phi\rangle dt=\langle\theta_0,\phi(0,x)\rangle+\int_0^T\langle S,\phi\rangle dt,
\end{align*}
where $u=u\Big|_{\nu=0}$. A weak solution $\theta(t)$ to \eqref{MG equation}-\eqref{MG Fourier symbol} on $[0,T]$ with $\nu=0$ is called a ``{\it vanishing viscosity}''{\it solution} if there exist sequences $\nu_n\rightarrow0$ and $\{\theta^{\nu_n}\}$ such that $\{\theta^{\nu_n}\}$ are smooth solutions to \eqref{MG equation}-\eqref{MG Fourier symbol} as given by Theorem~\ref{existence MG thm} and $\theta^{\nu_n}\rightarrow\theta$ in $C_w([0,T];L^2)$ as $\nu_n\rightarrow0$. 
\end{definition}

We prove that the system \eqref{MG equation}-\eqref{MG Fourier symbol} driven by a force $S$ possesses a compact global attractor in $L^2(\mathbb{T}^3)$ which is {\it upper semicontinuous} at $\nu=0$. More precisely, we have

\begin{thm}[Friedlander and Suen \cite{FS18}]\label{existence of MG attractor thm}
Assume $S\in C^\infty$. Then the system \eqref{MG equation}-\eqref{MG Fourier symbol} with $\nu=0$ possesses a compact global attractor $\mathcal{A}$ in $L^2(\mathbb{T}^3)$, namely
\begin{align*}
\mbox{$\mathcal{A}=\{\theta_0:\theta_0=\theta(0)$ for some bounded complete ``vanishing viscosity'' solution $\theta(t)\}$.}
\end{align*}
For any bounded set $\mathcal{B}\subset L^2(\mathbb{T}^3)$, and for any $\varepsilon,T>0$, there exists $t_0$ such that for any $t_1>t_0$, every ``vanishing viscosity'' solution $\theta(t)$ with $\theta(0)\in\mathcal{B}$ satisfies 
\begin{align*}
\|\theta(t)-x(t)\|_{L^2}<\varepsilon, \forall t\in[t_1,t_1+T],
\end{align*}
for some complete trajectory $x(t)$ on the global attractor $(x(t)\in\mathcal{A}, \forall t\in(-\infty,\infty))$. Furthermore, for $\nu\in[0,1]$, there exists a compact global attractor $\mathcal{A}^{\nu}\subset L^2$ for \eqref{MG equation}-\eqref{MG Fourier symbol} such that $\mathcal{A}^{0}=\mathcal{A}$ and $\mathcal{A}^{\nu}$ is {\it upper semicontinuous} at $\nu=0$, which means that
\begin{align}\label{uc_0}
\mbox{$\sup_{\phi\in\mathcal{A}^{\nu}}\inf_{\psi\in\mathcal{A}}\|\phi-\psi\|_{L^2}\rightarrow0$ as $\nu\rightarrow0$.}
\end{align}
\end{thm}

We give a brief discussion of the proof of Theorem~\ref{existence of MG attractor thm}, and we refer the interested reader to \cite{FS18} for full details. Roughly speaking, Theorem~\ref{existence of MG attractor thm} can be divided into two parts, namely:
\begin{itemize}
\item[1.] existence of global attractors $\mathcal{A}^{\nu}$ for \eqref{MG equation}-\eqref{MG Fourier symbol} with $\nu\ge0$; and 
\item[2.] upper semicontinuity of $\mathcal{A}^{\nu}$ at $\nu=0$.
\end{itemize}

For the existence of global attractors $\mathcal{A}^{\nu}$, it can be proved in the following several steps:

\begin{itemize}
\item[I.] First let $\theta(t)$ be a ``vanishing viscosity'' solution of \eqref{MG equation}-\eqref{MG Fourier symbol} on $[0,\infty)$ with $\theta(0)\in L^2$. Then $\theta(t)$ satisfies the following energy equality:
\begin{align}\label{energy equality MG}
\frac{1}{2}\|\theta(t)\|_{L^2}^2+\kappa\int_{t_0}^t\|\nabla\theta(s,\cdot)\|^2_{L^2}ds=\frac{1}{2}\|\theta(t_0)\|^2_{L^2}+\int_{t_0}^t\int_{\mathbb{T}^3}S\theta dxds,
\end{align}
for all $0\le t_0\le t$. The energy equality implies that:
\begin{itemize}
\item Every ``vanishing viscosity'' solution to \eqref{MG equation}-\eqref{MG Fourier symbol} is strongly continuous in time $t$.
\item There exists an absorbing ball $\mathcal{Y}$ for \eqref{MG equation}-\eqref{MG Fourier symbol} given by
\begin{align}\label{absorbing ball MG}
\mathcal{Y}=\{\theta\in L^2: \|\theta\|_{L^2}\le R\},
\end{align}
where $R$ is any number larger than $\kappa^{-1}\|S\|_{H^{-1}(\mathbb{T}^3)}$. 
\end{itemize}
\item[II.] Next we define $\pi^{\nu}: L^2\rightarrow L^2$ as the map $\pi^{\nu}\theta_0=\theta^{\nu}$, where $\theta^\nu$ is the solution to \eqref{MG equation}-\eqref{MG Fourier symbol} given by Theorem~\ref{existence MG thm}. Then using \cite[Lemma~6.7]{FS18}, for $t>0$, $\pi^{\nu}(t)\theta_0$ is continuous in $\nu$, uniformly for $\theta_0$ in compact subsets of $L^2$.
\item[III.] We denote the weak distance on $L^2(\mathbb{T}^3)$ by
\begin{align*}
d_w(\phi,\psi)=\sum_{k\in\mathbb{Z}^3}\frac{1}{2^{|k|}}\frac{|\hat\phi_k-\hat\psi_k|}{1+|\hat\phi_k-\hat\psi_k|},
\end{align*}
where $\hat\phi_k$ and $\hat\psi_k$ are the Fourier coefficients of $\phi$ and $\psi$. If 
\begin{align*}
\mbox{$\mathcal{E}[T,\infty)=\{\theta(\cdot):\theta(\cdot)$ is a ``vanishing viscosity'' solution of \eqref{MG equation}-\eqref{MG Fourier symbol}}\\
\mbox{on $[T,\infty)$ and $\theta\in\mathcal{Y}$ for all $t\in[T,\infty)\}$},
\end{align*}
\begin{align*}
\mbox{$\mathcal{E}(-\infty,\infty)=\{\theta(\cdot):\theta(\cdot)$ is a ``vanishing viscosity'' solution of \eqref{MG equation}-\eqref{MG Fourier symbol}}\\
\mbox{on $(-\infty,\infty)$ and $\theta\in\mathcal{Y}$ for all $t\in(-\infty,\infty)\}$},
\end{align*}
then $\mathcal{E}$ is an {\it evolutionary system} (see \cite{C09} and \cite{CD14} for the definition), so by \cite[Theorem~4.5]{CD14}, there exists a weak global attractor $\mathcal{A}_w$ to $\mathcal{E}$ with
\begin{align*}
\mbox{$\mathcal{A}_w = \{\theta_0 : \theta_0 = \theta(0)$ for some $\theta \in \mathcal{E}((-\infty,\infty))\}$}. 
\end{align*}
Applying the arguments given in \cite{CD14}, $\mathcal{E}$ satisfies all the following properties: 
\begin{itemize}
\item $\mathcal{E}([0, \infty))$ is a compact set in $C([0, \infty); \mathcal{Y}_w)$, here $\mathcal{Y}_w$ refers to the metric space $(\mathcal{Y}, d_w)$;
\item for any $\varepsilon>0$, there exists $\delta>0$ such that for every $\theta\in\mathcal{E}([0,\infty))$ and $t >0$, $$\|\theta(t)\|_{L^2}\le \|\theta(t)\|_{L^2} + \varepsilon,$$ for $t_0$ a.e. in $(t-\delta,t)\cap[0,\infty)$;
\item if $\theta_n\in\mathcal{E}([0,\infty))$ and $\theta_n\rightarrow\theta\in\mathcal{E}([0,\infty))$ in $C([0, \infty); \mathcal{Y}_w)$ for some $T>0$, then $\theta_n(t)\rightarrow\theta(t)$ strongly a.e. in $[0,T]$.
\end{itemize}
Together with \cite[Theorem~4.5]{CD14}, it implies that the strong global attractor $\mathcal{A}_s$ for \eqref{MG equation}-\eqref{MG Fourier symbol} with $\nu=0$ exists, it is strongly compact and $\mathcal{A}^0:=\mathcal{A}_s = \mathcal{A}_w$. The case for \eqref{MG equation}-\eqref{MG Fourier symbol} with $\nu>0$ is just similar.
\end{itemize}

On the other hand, to prove the upper semicontinuity of $\mathcal{A}^{\nu}$ at $\nu=0$, we note that the absorbing ball $\mathcal{Y}$ as given by \eqref{absorbing ball MG} has radius which is independent of $\nu$. Hence for all $\nu\ge0$, $\mathcal{A}^{\nu}$ satisfies
\begin{itemize}
\item $\pi^{0}(t)\mathcal{A^\nu}=\mathcal{A^\nu}$ for all $t\in\R$; 
\item for any bounded set $\mathcal{B}$, $\sup_{\phi\in\pi^{0}(t)\mathcal{B}}\inf_{\psi\in\mathcal{A}}d_{w}(\phi,\psi)\rightarrow0$ as $t\rightarrow0$.
\end{itemize}
We also have that $\mathcal{E}([0, \infty))$ is a compact set in $C([0, \infty); \mathcal{Y}_w)$ such that $\mathcal{A}^{\nu}\subset \mathcal{K}$ for every $\nu\in[0,1]$. Together with the continuity of $\pi^{\nu}(t)\theta_0$ in $\nu$, the result from \cite{HOR15} implies the {\it weak} upper semicontinuity, namely
\begin{align}\label{weak upper continuity MG diffusive}
\mbox{$\sup_{\phi\in\mathcal{A}^\nu}\inf_{\psi\in\mathcal{A}}d_w(\phi,\psi)\rightarrow0$ as $\nu\rightarrow0$.}
\end{align}
Moreover, for any $\phi^{\nu_j}\in\mathcal{A}^{\nu_j}$ and $\psi_j\in\mathcal{A}$, if
\begin{align*}
\lim_{j\rightarrow\infty} d_w(\phi^{\nu_j},\psi_j)=0,
\end{align*}
for some sequence $\nu_j\rightarrow0$, then
\begin{align*}
\lim_{j\rightarrow0}\|\phi^{\nu_j}-\psi_j\|_{L^2}=0.
\end{align*}
In other words, the {\it weak} upper semicontinuity implies the {\it strong} upper semicontinuity, and hence \eqref{weak upper continuity MG diffusive} further implies the {\it strong} upper semicontinuity of $\mathcal{A}^{\nu}$ at $\nu=0$. 

\section{Non-diffusive active scalar equations}\label{s:Non-diffusive active scalar equations}

We switch our attention to the non-diffusive active scalar equations \eqref{abstract active scalar eqn}-\eqref{zero mean assumption} for $\kappa=0$ and $\nu\ge0$, namely
\begin{align}
\label{abstract active scalar eqn non-diffusive} \left\{ \begin{array}{l}
\partial_t\theta+u\cdot\nabla\theta=S, \\
u_j[\theta]=\partial_{x_i} T_{ij}^{\nu}[\theta],\theta(x,0)=\theta_0(x)
\end{array}\right.
\end{align}
with $\nu\ge0$, and $S$, $\theta_0$ satisfy the zero mean assumptions \eqref{zero mean assumption on data and forcing}-\eqref{zero mean assumption}. We address the well-posedness of \eqref{abstract active scalar eqn non-diffusive} and convergence of solutions as $\nu\to0$. The results can then be applied to the thermally non diffusive MG$^\nu$ equation \eqref{MG equation}-\eqref{MG Fourier symbol} (with $\kappa=0$) as well as the IPMB$^\nu$ equation \eqref{IPMB constitutive law}-\eqref{Fourier multiplier symbol IPMB}.

\subsection{Well-posedness in various spaces for $\nu>0$}\label{s:Well-posedness in various spaces nu+ non-diffusive}

In \cite{FV11b}, Friedlander and Vicol showed that the singular inviscid MG$^0$ equation is {\it ill-possed} in the sense of Hadamard in any Sobolev space. One of the main reasons for the ill-posedness is that, the Fourier multiplier symbols $\widehat{M}^0$ given by \eqref{MG Fourier symbol} is an {\it even} function in $k$. When we perform energy estimates on \eqref{MG equation}, we may have trouble on controlling the following term
\[R:=\int(-\Delta)^\frac{s}{2}u\cdot\nabla\theta(-\Delta)^\frac{s}{2}\theta.\]
Without the diffusive term $\kappa \Delta \theta$, we loss control on the term $R$ and the only hope to treat it would be to discover a commutator structure. If the operator $\widehat{\partial_{x_i} T_{ij}^{\nu}}(k)$ is odd, then there is an extra {\it cancellation} which allows us to close the estimates (at the level of Sobolev spaces) by the Coifman-Meyer type commutator estimate. On the other hand, for the case when the operator $\widehat{\partial_{x_i} T_{ij}^{\nu}}(k)$ is even, such argument fails and one has to seek for some other type of estimates. In the case of $\nu>0$, however, the operators $\partial_x T^\nu$ are smoothing order 2 which give rise to well-posedness for \eqref{abstract active scalar eqn} with $\theta_0\in W^{s,d}$ for $s\ge0$ and smooth forcing term $S$. The results are summarised as follows:

\begin{thm}[Friedlander and Suen \cite{FS19}]\label{Well-posedness in Sobolev space inviscid thm}
Let $\theta_0\in W^{s,d}$ for $s\ge0$ and $S$ be a $C^\infty$-smooth source term. Then for each $\nu>0$, we have:
\begin{itemize}
\item if $s=0$, there exists unique global weak solution to \eqref{abstract active scalar eqn non-diffusive} such that
\begin{align*}
\theta^\nu&\in BC((0,\infty);L^d),\\
u^\nu&\in C((0,\infty);W^{2,d}).
\end{align*}
In particular, $\theta^\nu(\cdot,t)\rightarrow\theta_0$ weakly in $L^d$ as $t\rightarrow0^+$. Here $BC$ stands for {\it bounded continuous functions}.
\item if $s>0$, there exists a unique global-in-time solution $\theta^\nu$ to \eqref{abstract active scalar eqn non-diffusive} such that $\theta^\nu(\cdot,t)\in W^{s,d}$ for all $t\ge0$. Furthermore, for $s=1$, we have the following single exponential growth in time on $\|\nabla\theta^\nu(\cdot,t)\|_{L^d}$:
\begin{align*}
\|\nabla\theta^\nu(\cdot,t)\|_{L^d}\le C \|\nabla\theta_0\|_{L^d}\exp\left(C\left(t\|\theta_0\|_{W^{1,d}}+t^2\|S\|_{L^\infty}+t\|S\|_{W^{1,d}}\right)\right),
\end{align*}
where $C>0$ is a constant which depend only on $\nu$ and the spatial dimension $d$.
\end{itemize}
\end{thm}

The proof of Theorem~\ref{Well-posedness in Sobolev space inviscid thm} for the case of $s=0$ relies on the existence and uniqueness of the flow map, which is essential for Euler system as well \cite{BK12}. We briefly sketch here and the full details can be found in \cite{FS19}.

In view of condition (A3), by applying Fourier multiplier theorem (see Stein \cite{S70}), given $p>1$, there exists some constant $C=C(\nu,p,d)>0$ such that
\begin{align}\label{Fourier multiplier theorem inviscid}
\|u^\nu(\cdot,t)\|_{W^{2,p}}\le C \|\theta^\nu(\cdot,t)\|_{L^p}.
\end{align}
Together with \eqref{Fourier multiplier theorem inviscid} and embedding theorems, one can show that the {\it Log-Lipschitzian} norm of $u^\nu$ given by $\|u^\nu(\cdot,t)\|_{L.L.}$ is bounded in terms of $\theta_0$ and $S$:
\begin{align}\label{LL bound}
\|u^\nu(\cdot,t)\|_{L.L.}\le C\left(\|\theta_0\|_{L^d}+t\|S\|_{L^\infty}\right).
\end{align} 
Next, we consider the standard mollifier $\rho\in C^\infty_0 $, and we set $\theta_{{(n)},0}=\rho_n *\theta_0$ for $n\in\N$ and $\rho_n(x)=n^d\rho(nx)$. By a standard argument, given $\nu>0$, we can obtain a sequence of global smooth solution $(\theta_{(n)}^\nu,u_{(n)}^\nu)$ to \eqref{abstract active scalar eqn non-diffusive} with $\theta_{(n)}^\nu(x,0)=\theta_{{(n)},0}$ and $u_{(n)}^\nu=\partial_{x_i} T_{ij}^{\nu}[\theta_{(n)}^\nu]$. We define $\psi_n(x,t)$ to be the flow map given by
\begin{align*}
\partial_t\psi_n(x,t)=u_{(n)}^\nu(\psi_n(x,t),t),
\end{align*}
then $\psi_n$ satisfies
\begin{align}\label{bound on flow map}
\|\psi_n(\cdot,t)\|_{*}\le C \exp\left(\int_0^t\|u_{(n)}^\nu(\cdot,\tilde{t})\|_{L.L.}d\tilde{t}\right),
\end{align}
where the norm $\|\cdot\|_{*}$ is given by
\begin{align*}
\|\psi\|_{*}=\sup_{x\neq y}\Phi(|\psi(x)-\psi(y)|,|x-y|)
\end{align*}
with
\begin{align*}
\Phi(r,s)=\left\{ \begin{array}{l}
\mbox{$\max\{\frac{1+|\log(s)|}{1+|\log(r)|},\frac{1+|\log(r)|}{1+|\log(s)|}\}$, if $(1-s)(1-r)\ge0$,}\\
\mbox{$(1+|\log(s)|)(1+|\log(r)|)$, if $(1-s)(1-r)\le0$.}
\end{array}\right.
\end{align*}
Using \eqref{LL bound} and \eqref{bound on flow map} (with $u^\nu$ replaced by $u_{(n)}^\nu$), we obtain
\begin{align}\label{bound on psi 1}
|\psi_n(x_1,t)-\psi_n(x_2,t)|\le \alpha(t)|x_1-x_2|^{\beta(t)}
\end{align}
for all $(x_1,t),(x_2,t)\in\R^d\times\R^+$, where $\alpha(t),\beta(t)$ are some continuous functions which depend on $\theta_0$ and $S$. Furthermore, for $t_1,t_2\ge0$, using \eqref{Fourier multiplier theorem inviscid} (with $u^\nu$ replaced by $u_{(n)}^\nu$),
\begin{align}\label{bound on psi 2}
|\psi_n(x,t_1)-\psi_n(x,t_2)|\le C |t_2-t_1|(\|\theta_0\|_{L^p}+\max\{t_1,t_2\}\|S\|_{L^\infty}).
\end{align}
The estimates \eqref{bound on psi 1} and \eqref{bound on psi 2} imply that the family $\{\psi_n\}_{n\in\N}$ is bounded and equicontinuous on every compact set in $\R^d\times\R^+$. By the Arzela-Ascoli theorem, there exists a limiting trajectory $\psi(x,t)$ as $n\rightarrow\infty$. Performing the same analysis for $\{\psi_n^{-1}\}$, where $\psi_n^{-1}$ is the inverse of $\psi_n$, we see that $\psi(x,t)$ is a Lebesgue measure preserving homeomorphism as well. Define $\theta^\nu(x,t)=\theta_0(\psi^{-1}(x,t))$ and $u^\nu=\partial_{x_i} T_{ij}^{\nu}[\theta^\nu]$. Then one can show that $(\theta^\nu,u^\nu)$ is a weak solution to \eqref{abstract active scalar eqn}. To show that $(\theta^\nu,u^\nu)$ is unique, let $T>0$ and $\nu>0$, and suppose that $(\theta^{\nu,1},u^{\nu,1})$ and $(\theta^{\nu,2},u^{\nu,2})$ solve \eqref{abstract active scalar eqn} on $\mathbb{T}^d\times[0,T]$ with $\theta^{\nu,1}(\cdot,0)=\theta^{\nu,2}(\cdot,0)=\theta_0$. Then there exists a constant $C>0$ such that for all $\delta\in(0,1)$ and $k\in\{-1\}\cup\N$, we have
\begin{align*}
&\|\Delta_k (\theta^{\nu,1}-\theta^{\nu,2})(\cdot,t)\|_{L^\infty}\\
&\le 2^{k\delta}(k+1)C(\|u^{\nu,1}(\cdot,t)\|_{\overline{L.L.}}+\|u^{\nu,2}(\cdot,t)\|_{\overline{L.L.}})\|(\theta^{\nu,1}-\theta^{\nu,2})(\cdot,t)\|_{B^{-\delta}_{\infty,\infty}}, \forall t\in[0,T],
\end{align*}
where $\|\cdot\|_{\overline{L.L.}}=\|\cdot\|_{L^\infty}+\|\cdot\|_{L.L.}$ and $\Delta_k$'s are the dyadic blocks for $k\in\{-1\}\cup\N$. We define
\begin{align*}
\bar{t}=\sup\left\{t\in[0,T]:C\int_0^t(\|u^{\nu,1}(\cdot,\tilde{t})\|_{\overline{L.L.}}+\|u^{\nu,2}(\cdot,\tilde{t})\|_{\overline{L.L.}})d\tilde{t}\le\frac{1}{2}\right\},
\end{align*}
then by the bounds \eqref{Fourier multiplier theorem inviscid} and \eqref{LL bound}, $\bar{t}$ is well-defined. We let 
\begin{align*}
\delta_{\bar{t}}=C\int_0^{\bar{t}}(\|u^{\nu,1}(\cdot,\tilde{t})\|_{\overline{L.L.}}+\|u^{\nu,2}(\cdot,\tilde{t})\|_{\overline{L.L.}})d\tilde{t}.
\end{align*} 
Using \cite[Theorem~3.28]{BCD11}, for all $k\ge-1$ and $t\in[0,\bar{t}]$,
\begin{align*}
2^{-k\delta_{\bar{t}}}\|\Delta_k (\theta^{\nu,1}-\theta^{\nu,2})(\cdot,t)\|_{L^\infty}\le \frac{1}{2}\sup_{t\in[0,\bar{t}]}\|(\theta^{\nu,1}-\theta^{\nu,2})(\cdot,t)\|_{B^{-\delta_{\bar{t}}}_{\infty,\infty}}.
\end{align*}
Summing over $k$ and taking supremum over $[0,\bar{t}]$, we conclude that $\theta^{\nu,1}=\theta^{\nu,2}$ on $[0,\bar{t}]$. By repeating the argument a finite number of times, we obtain the uniqueness on the whole interval $[0,T]$. This concludes our sketch of the proof of Theorem~\ref{Well-posedness in Sobolev space inviscid thm}.

Next we study the Gevrey-class $s$ solutions to \eqref{abstract active scalar eqn non-diffusive} for $\nu>0$ when the initial datum $\theta_0$ and forcing term $S$ are in the same Gevrey-class. We prove that there exists a unique global-in-time Gevrey-class $s$ solution $\theta^\nu$ with radius of convergence bounded below by some positive function $\tau(t)$ for all $t\in[0,\infty)$. More precisely, we have:

\begin{thm}[Friedlander and Suen \cite{FS19}]\label{Gevrey-class global well-posedness inviscid thm}
Fix $s\ge1$. Let $\theta_0$ and $S$ be of Gevrey-class $s$ with radius of convergence $\tau_0>0$. Then there exists a unique Gevrey-class $s$ solution $\theta^\nu$ to \eqref{abstract active scalar eqn non-diffusive} on $\mathbb{T}^d\times[0,\infty)$ with radius of convergence at least $\tau=\tau(t)$ for all $t\in[0,\infty)$, where $\tau$ is a decreasing function satisfying
\begin{align}\label{lower bound on tau}
\tau(t)\ge\tau_0e^{-C\left(\|e^{\tau_0\Lambda^\frac{1}{s}}\theta_0\|_{L^2}+2\|e^{\tau_0\Lambda^\frac{1}{s}}S\|_{L^2}\right)t}.
\end{align}
Here $C>0$ is a constant which depends on $\nu$ but independent of $t$.
\end{thm}

The Gevrey-class $s$ is given by
\begin{align*}
\bigcup_{\tau>0}\mathcal{D}(\Lambda^re^{\tau\Lambda^\frac{1}{s}}),
\end{align*}
for any $r\ge0$, where
\begin{align*}
\|\Lambda^re^{\tau\Lambda^\frac{1}{s}}f\|^2_{L^2}=\sum_{k\in\mathbb{Z}^d_*}|k|^{2r}e^{2\tau|k|^\frac{1}{s}}|\hat{f}(k)|^2,
\end{align*}
where $\tau=\tau(t)>0$ denotes the radius of convergence and $\Lambda=(-\Delta)^\frac{1}{2}$. Hence by taking $L^2$-inner product of \eqref{abstract active scalar eqn non-diffusive}$_1$ with $e^{2\tau\Lambda^{\frac{1}{s}}}\theta^\nu$ and applying H\"{o}lder's inequality, we obtain
\begin{align*}\label{a priori estimates on theta}
\frac{1}{2}\frac{d}{dt}\|e^{\tau\Lambda^{\frac{1}{s}}}\theta^\nu\|^2_{L^2}-\dot{\tau}\|\Lambda^\frac{1}{2s}e^{\tau\Lambda^{\frac{1}{s}}}\theta^\nu\|^2_{L^2}&\le\Big|-\langle u^\nu\cdot\nabla\theta^\nu,e^{2\tau\Lambda^{\frac{1}{s}}}\theta^\nu\rangle\Big|+\|e^{\tau\Lambda^\frac{1}{s}}S\|_{L^2}\|e^{\tau\Lambda^\frac{1}{s}}\theta^\nu\|_{L^2}\notag.
\end{align*}
The key step for proving Theorem~\ref{Gevrey-class global well-posedness inviscid thm} is to estimate the term $\Big|-\langle u^\nu\cdot\nabla\theta^\nu,e^{2\tau\Lambda^{\frac{1}{s}}}\theta^\nu\rangle\Big|$. Using a Cauchy-Kowalewski-type argument and together with condition (A3), it can be showed in \cite[Lemma~4.1]{FS19} that
\begin{equation*}
\left|-\langle e^{\tau\Lambda^{\frac{1}{s}}}(u^\nu\cdot\nabla\theta^\nu),e^{\tau\Lambda^{\frac{1}{s}}}\theta^\nu\rangle\right|\le C\tau\|e^{\tau\Lambda^\frac{1}{s}}\theta^\nu\|_{L^2}\|\Lambda^\frac{1}{2s}e^{\tau\Lambda^\frac{1}{s}}\theta^\nu\|^2_{L^2},
\end{equation*}
and we obtain that
\begin{align*}
\frac{1}{2}\frac{d}{dt}\|e^{\tau\Lambda^{\frac{1}{s}}}\theta^\nu\|^2_{L^2}-\dot{\tau}\|\Lambda^{\frac{1}{2s}}e^{\tau\Lambda^{\frac{1}{s}}}\theta^\nu\|^2_{L^2}+\kappa\|e^{\tau\Lambda^\frac{1}{s}}\theta^{\nu}\|^2_{L^2}\le C\tau\|e^{\tau\Lambda^\frac{1}{s}}\theta^\nu\|_{L^2}\|\Lambda^\frac{1}{2s}e^{\tau\Lambda^\frac{1}{s}}\theta^\nu\|^2_{L^2}.
\end{align*}
By choosing $\tau>0$ such that
\begin{equation*}
\dot{\tau}+C\tau\|e^{\tau\Lambda^\frac{1}{s}}\theta^\nu\|_{L^2}=0,
\end{equation*}
we have
\begin{equation*}
\|e^{\tau(t)\Lambda^\frac{1}{s}}\theta^\nu(t)\|_{L^2}\le\|e^{\tau_0\Lambda^\frac{1}{s}}\theta_0\|_{L^2}+2\|e^{\tau_0\Lambda^\frac{1}{s}}S\|_{L^2}.
\end{equation*}
and $\tau$ satisfies the lower bound \eqref{lower bound on tau}. We remark that, for the diffusive case given by the system \eqref{abstract active scalar eqn diffusive}, one can obtain global-in-time existence of solution Gevrey class $s\ge1$ with lower bound on $\tau(t)$ that does not vanish as $t\rightarrow\infty$, refer to \cite[Remark~4.3]{FS19} for more details. 

\subsection{Well-posedness in various spaces for $\nu=0$}\label{s:Well-posedness in various spaces nu0 non-diffusive}

In this subsection we study the non-diffusive equations \eqref{abstract active scalar eqn non-diffusive} for $\nu=0$:
\begin{align}
\label{abstract active scalar eqn nu=0 non-diffusive} \left\{ \begin{array}{l}
\partial_t\theta^0+u^0\cdot\nabla\theta^0=S, \\
u^0_j=\partial_{x_i} T^0_{ij}[\theta^0],\theta^0(x,0)=\theta_0(x).
\end{array}\right.
\end{align}
When $\nu = 0$ and condition (A2) is imposed, as it was proved in \cite{FV11b}, the equation \eqref{abstract active scalar eqn nu=0 non-diffusive} is {\it ill-posed} in the sense of Hadamard, which means that the solution map associated to the Cauchy problem for \eqref{abstract active scalar eqn nu=0 non-diffusive} is not Lipschitz continuous with respect to perturbations in the initial datum around a specific steady profile $\theta_0$, in the topology of a certain Sobolev space $X$. Nevertheless, as pointed out in \cite{FV11b}, it is possible to obtain the local existence and uniqueness of solutions to \eqref{abstract active scalar eqn nu=0 non-diffusive} in spaces of real-analytic functions, owing to the fact that the derivative loss in the nonlinearity $u^0\cdot\nabla\theta^0$ is of order at most one (both in $u^0$ and in $\nabla\theta^0$). In \cite{FS19}, we extended the results of \cite{FV11b} to the case of Gevrey-class solutions which are summarised as follows:

\begin{thm}[Friedlander and Suen \cite{FS19}]\label{Gevrey-class local well-posedness nu0 inviscid thm}
Fix $s\ge1$, $r>\frac{d}{2}+\frac{3}{2}$ and $K_0>0$. Let $\theta^0(\cdot,0)=\theta_0$ and $S$ be of Gevrey-class $s$ with radius of convergence $\tau_0>0$ and satisfy
\begin{align}\label{bounds on gevrey norm of S and theta0'}
\|\Lambda^re^{\tau_0\Lambda^\frac{1}{s}}\theta^0(\cdot,0)\|_{L^2}\le K_0,\qquad\|\Lambda^re^{\tau_0\Lambda^\frac{1}{s}}S\|_{L^2}\le K_0.
\end{align}
For $\nu=0$, under the condition {\rm (A2)}, there exists $\bar{T},\bar{\tau}>0$ and a unique Gevrey-class $s$ solution $\theta^0$ to \eqref{abstract active scalar eqn nu=0 non-diffusive} defined on $\mathbb{T}^d\times[0,\bar{T}]$ with radius of convergence at least $\bar{\tau}$. In particular, there exists a constant $C=C(K_0)>0$ such that for all $t\in[0,\bar{T}]$,
\begin{align}
\|\Lambda^re^{\bar{\tau}\Lambda^\frac{1}{s}}\theta^0(\cdot,t)\|_{L^2}\le C.\label{bound on solution 2}
\end{align}
The bound \eqref{bound on solution 2} also applies on $\theta^\nu$ for $\nu>0$.
\end{thm}

In contrast, when $\nu = 0$ and condition (A2$^*$) is in force, the operator $\partial_x T^0$ becomes a zero order operator with $\partial_x T^0:L^2\rightarrow L^2$ being bounded. Following the idea given in \cite{FRV12}, we show that the equation \eqref{abstract active scalar eqn nu=0 non-diffusive} is locally well-posed in Sobolev space $H^s$ for $s>\frac{d}{2}+1$: 

\begin{thm}[Friedlander and Suen \cite{FS19}]\label{Local well-posedness in Sobolev space nu0 inviscid thm}
For $d\ge2$, we fix $s>\frac{d}{2}+1$. Assume that $\theta_0,S\in H^s(\mathbb{T}^d)$ have zero-mean on $\mathbb{T}^d$. Then for $\nu=0$, under the condition {\rm (A2$^*$)}, there exists a $T>0$ and a unique smooth solution $\theta^0$ to \eqref{abstract active scalar eqn nu=0 non-diffusive} such that
$$\theta^0\in L^\infty(0,T;H^s(\mathbb{T}^d)).$$
\end{thm}

The proof of Theorem~\ref{Local well-posedness in Sobolev space nu0 inviscid thm} consists of three steps, which can be briefly outline as follows (details can be found in \cite{FS19}):
\begin{itemize}
\item[I.] We first construct a sequence of approximations $\{\theta_n\}_{n\ge1}$ given by the solutions of
\begin{align*}
\partial_t\theta_1&=S\notag\\
\theta_1(\cdot,0)&=\theta_0.
\end{align*}
and for $n>1$,
\begin{align}\label{theta n=n}
\partial_t\theta_n+u_{n-1}\cdot\nabla\theta_n&=S\notag\\
\theta_{n-1}&=\partial_xT^0[\theta_{n-1}]\\
\theta_n(\cdot,0)&=\theta_0.\notag
\end{align}
Then by applying \cite[Theorem~A1]{FRV12}, one can show that $\theta_n\in L^\infty(0,T;H^s)$ for all $n\in\mathbb{N}$.
\item[II.] Next, by induction on $n$, we prove that $\|\Lambda^s \theta_{n}(\cdot,t)\|_{L^2}$ is bounded on $[0,T]$ for some $T>0$. Assume that
\begin{align*}
\|\Lambda^s\theta_j\|_{L^\infty(0,T;L^2)}\le C\|\Lambda^s\theta_0\|_{L^2},
\end{align*}
for $1\le j\le n-1$. We apply $\Lambda^s$ on \eqref{theta n=n} and take inner product with $\Lambda^s\theta_n$ to obtain
\begin{align}\label{estimate on theta n identity}
\frac{1}{2}\frac{d}{dt}\intox|\Lambda^s\theta_n|^2+\intox\Lambda^s\theta_n\cdot\Lambda^s(u_{n-1}\cdot\nabla\theta_n)=\intox\Lambda^s\theta_n\cdot\Lambda^s S.
\end{align}
Using commutator estimates (see \cite{FRV12} for example), the term involving $u_{n-1}$ can be bounded in terms of $\theta_n$ and $\theta_{n-1}$:
\begin{align*}
\left|\intox\Lambda^s\theta_n\cdot\Lambda^s(u_{n-1}\cdot\nabla\theta_n)\right|\le C\|\Lambda^s\theta_n\|_{L^2}(\|\Lambda^s \theta_{n-1}\|_{L^2}\|\Lambda^s\theta_n\|_{L^2}).
\end{align*}
Hence by integrating the identity \eqref{estimate on theta n identity} over $t$, choosing $T$ small enough and applying the induction hypothesis, $\|\Lambda^s \theta_{n}(\cdot,t)\|_{L^2}$ is bounded on $[0,T]$ as well.
\item[III.] Finally, we show that $\{\theta_n\}_{n\ge0}$ is a Cauchy sequence. This can be done by considering the difference $\tilde\theta_{n}=\theta_n-\theta_{n-1}$ and one can prove that
\begin{align*}
\frac{d}{dt}\|\Lambda^{s-1}\tilde\theta_n\|_{L^2}\le C(\|\Lambda^s\theta_0\|_{L^2}\|\Lambda^{s-1}\tilde\theta_n\|_{L^2}+\|\Lambda^{s-1}\tilde\theta_{n-1}\|_{L^2}\|\Lambda^s\theta_0\|_{L^2}).
\end{align*}
Integrating the above over $t$ and choosing $T$ small enough, we obtain
\begin{align*}
\sup_{t\in[0,T]}\|\Lambda^{s-1}\tilde\theta_{n}(\cdot,t)\|_{L^2}\le\frac{1}{2}\sup_{t\in[0,T]}\|\Lambda^{s-1}\tilde\theta_{n-1}(\cdot,t)\|_{L^2}.
\end{align*}
Thus $\theta_n$ is Cauchy in $L^\infty(0,T;H^{s-1})$ with $\theta_n$ converges strongly to $\theta^0$ in $L^\infty(0,T,H^{s-1})$. Since we assume that $s>\frac{d}{2}+1$, this also implies that the strong convergence occurs in a H\"{o}lder space relative to $x$ as $n\rightarrow\infty$, hence the limiting function $\theta^0$ is a solution of \eqref{abstract active scalar eqn nu=0 non-diffusive}. Uniqueness of $\theta^0$ follows by the same argument given in \cite{FRV12} and we omit the details.
\end{itemize}

\subsection{Convergence of solutions as $\nu\to0$}\label{s:convergence of solutions inviscid}

In this subsection, we address the convergence of solutions to \eqref{abstract active scalar eqn non-diffusive} as $\nu\rightarrow0$. Depending on the conditions (A2) and (A2$^*$), we can address the convergence of solutions in two cases respectively:

\subsubsection{Gevrey-class solutions:} We focus on the case for Gevrey-class solutions $\theta^\nu$ to \eqref{abstract active scalar eqn non-diffusive} when (A2) is in force. By Theorem~\ref{Gevrey-class local well-posedness nu0 inviscid thm}, given Gevrey-class $s$ initial datum $\theta_0$ and forcing $S$, there exists $\bar{T},\bar{\tau}>0$ and a unique Gevrey-class solution $\theta^\nu$ to \eqref{abstract active scalar eqn non-diffusive} defined on $[0,\bar{T}]$ with radius of convergence at least $\bar{\tau}$ for all $\nu\ge0$. In particular, the Gevrey-class solutions $\theta^\nu$ converges to $\theta^0$ in some Gevrey-class norm as $\nu\rightarrow0$ and the results are summarised in the following theorem:

\begin{thm}[Friedlander and Suen \cite{FS19}]\label{Convergence of Gevrey solutions as nu goes to 0 inviscid thm} 
Under the condition {\rm (A2)}, if $\theta^\nu$ and $\theta^0$ are Gevrey-class $s$ solutions to \eqref{abstract active scalar eqn non-diffusive} for $\nu>0$ and $\nu=0$ respectively with initial datum $\theta_0$ on $\mathbb{T}^d\times[0,\bar{T}]$ with radius of convergence at least $\bar{\tau}$ as described in Theorem~\ref{Gevrey-class local well-posedness nu0 inviscid thm}, then there exists $T<\bar{T}$ and $\tau=\tau(t)<\bar{\tau}$ such that, for $t\in[0,T]$, we have
\begin{align}\label{convergence Gevrey}
\lim_{\nu\rightarrow0}\|(\Lambda^re^{\tau\Lambda^\frac{1}{s}}\theta^{\nu}-\Lambda^re^{\tau\Lambda^\frac{1}{s}}\theta^0)(\cdot,t)\|_{L^2}=0.
\end{align}
\end{thm}

The proof of Theorem~\ref{Convergence of Gevrey solutions as nu goes to 0 inviscid thm} relies on the estimates of the difference $\phi^\nu:=\theta^\nu-\theta^0$, and it can be shown that $\phi^\nu$ satisfies 
\begin{align*}
\frac{1}{2}\frac{d}{dt}\|\phi^\nu\|^2_{\tau,r}=\dot{\tau}\|\Lambda^\frac{1}{2s}\phi^\nu\|^2_{\tau,r}+\mathcal{R}_1+\mathcal{R}_2,
\end{align*}
where the terms $\mathcal{R}_1$ and $\mathcal{R}_2$ can be bounded as follows:
\begin{align*}
\mathcal{R}_1&\le C\|\Lambda^\frac{1}{2s}\theta^0\|_{\tau,r}\|\Lambda^\frac{1}{2s}\phi^\nu\|_{\tau,r}\|\phi^\nu\|_{\tau,r}\notag\\
&\qquad+C\|\Lambda^\frac{1}{2s}\theta^0\|_{\tau,r}\|\Lambda^\frac{1}{2s}\phi^\nu\|_{\tau,r}\left(\sum_{j\in\mathbb{Z}_*^d}|j|^{d+3}|\widehat{\theta^0}(j)|^2e^{2\tau|j|}|(\widehat{T^\nu}-\widehat{T^0})(j)|^2\right)^\frac{1}{2},\\
\mathcal{R}_2&\le C\|\Lambda^\frac{1}{2s}\phi^\nu\|_{\tau,r}^2\|\theta^\nu\|_{\tau,r}.
\end{align*}
By choosing $\tau=\tau(t)\le\bar{\tau}$ such that
\begin{align*}
\left\{ \begin{array}{l}
\dot{\tau}+C\|\theta^\nu\|_{\tau,r}+C\|\Lambda^\frac{1}{2s}\theta^0\|_{\tau,r}^2<0, \\
\tau<\bar{\tau},
\end{array}\right.
\end{align*}
and applying the bound \eqref{bound on solution 2} to conclude that
\begin{align*}
\frac{d}{dt}\|\phi^\nu\|_{\tau,r}^2\le C\|\phi^\nu\|_{\tau,r}^2+C\sum_{j\in\mathbb{Z}_*^d}|j|^{d+3}|\widehat{\theta^0}(j)|^2e^{2\tau|j|}|(\widehat{T^\nu}-\widehat{T^0})(j)|^2.
\end{align*}
Integrating the above with respect to $t$ and using the condition (A2), we have $\dis\lim_{\nu\to0}\|\phi^\nu\|_{\tau,r}=0$ and \eqref{convergence Gevrey} follows.

\subsubsection{$H^s$ solutions:} When condition (A2$^*$) is in force, by Theorem~\ref{Local well-posedness in Sobolev space nu0 inviscid thm}, the equation \eqref{abstract active scalar eqn non-diffusive} for $\nu=0$ is locally well-posed in Sobolev space $H^{s}$ for $s>\frac{d}{2}+1$. For sufficiently smooth initial data $\theta_0$ and forcing term $S$, one can show that $\|(\theta^{\nu}-\theta^0)(\cdot,t)\|_{H^s}\rightarrow0$ as $\nu\rightarrow0$ for $s>\frac{d}{2}+1$ and $t\in[0,T]$. Such result is parallel to the one proved in \cite{FS18}, in which the authors proved that if $\theta^{\nu},\theta^0$ are $C^\infty$ smooth classical solutions of the diffusive system \eqref{abstract active scalar eqn diffusive} for $\nu>0$ and $\nu=0$ respectively with initial datum $\theta_0\in L^2$ and forcing term $S\in C^\infty$, then $\|(\theta^{\nu}-\theta^0)(\cdot,t)\|_{H^s}\rightarrow0$ as $\nu\rightarrow0$ for $s\ge0$ and $t>0$. The convergence results are summarised below:

\begin{thm}[Friedlander and Suen \cite{FS19}]\label{Convergence of Sobolev solutions as nu goes to 0 inviscid thm} 
Under the condition {\rm (A2$^*$)}, we have
\begin{align}\label{convergence L2}
\lim_{\nu\rightarrow0}\|(\theta^{\nu}-\theta^0)(\cdot,t)\|_{H^{s-1}}=0,
\end{align}
and for $d\ge2$ and $s>\frac{d}{2}+1$ and $t\in[0,T]$, we have
\begin{align}\label{convergence sobolev}
\lim_{\nu\rightarrow0}\|(\theta^{\nu}-\theta^0)(\cdot,t)\|_{H^{s-1}}=0.
\end{align}
\end{thm}

It suffices to consider the case for the convergence in $L^2$ given by \eqref{convergence L2}, since the case for \eqref{convergence sobolev} follows by Gagliardo-Nirenberg interpolation inequality \cite{G59} and \cite{L59}. The key step of the proof is to estimate $\|(u^{\nu}-u^0)(\cdot,t)\|_{L^2}$, which can be bounded by $\|\phi^\nu(\cdot,t)\|_{L^2}^2+I(\nu,t)$ with $I(\nu,t)$ becoming zero as $\nu$ vanishes, see \cite{FS19} for further details.

\subsection{Applications to physical models}\label{s:applications to physical models inviscid}

We now apply our results discussed previous subsections to some physical models, namely the thermally non diffusive magnetogeostrophic (MG$^\nu$) equations \eqref{MG equation}-\eqref{MG Fourier symbol} with $\kappa=0$ and the incompressible porous media Brinkman equations (IPMB$^\nu$) \eqref{IPMB constitutive law}-\eqref{Fourier multiplier symbol IPMB}. The results are summarised in the following theorems (also refer to \cite{FS19} for details):

\begin{thm}[Well-posedness in Sobolev space for the MG$^\nu$ equations]\label{Well-posedness in Sobolev space MG non-viscid thm}
Let $\theta_0\in W^{s,3}$ for $s\ge0$ and $S$ be a $C^\infty$-smooth source term. Then for each $\nu>0$, we have:
\begin{itemize}
\item if $s=0$, there exists unique global weak solution to \eqref{MG equation}-\eqref{MG Fourier symbol} with $\kappa=0$ such that
\begin{align*}
\theta^\nu&\in BC((0,\infty);L^3),\\
u^\nu&\in C((0,\infty);W^{2,3}).
\end{align*}
In particular, $\theta^\nu(\cdot,t)\rightarrow\theta_0$ weakly in $L^3$ as $t\rightarrow0^+$.
\item if $s>0$, there exists a unique global-in-time solution $\theta^\nu$ to \eqref{MG equation}-\eqref{MG Fourier symbol} with $\kappa=0$ such that $\theta^\nu(\cdot,t)\in W^{s,3}$ for all $t\ge0$. Furthermore, for $s=1$, we have the following single exponential growth in time on $\|\nabla\theta^\nu(\cdot,t)\|_{L^3}$:
\begin{align*}
\|\nabla\theta^\nu(\cdot,t)\|_{L^3}\le C \|\nabla\theta_0\|_{L^3}\exp\left(C\left(t\|\theta_0\|_{W^{1,3}}+t^2\|S\|_{L^\infty}+t\|S\|_{W^{1,3}}\right)\right),
\end{align*}
where $C>0$ is a constant which depends only on some dimensional constants.
\end{itemize}
\end{thm}

\begin{thm}[Gevrey-class well-posedness for the MG$^\nu$ equations]\label{Gevrey-class global well-posedness MG non-viscid thm}
Fix $s\ge1$. Let $\theta_0$ and $S$ be of Gevrey-class $s$ with radius of convergence $\tau_0>0$. Then for each $\nu>0$, there exists a unique Gevrey-class $s$ solution $\theta^\nu$ to \eqref{MG equation}-\eqref{MG Fourier symbol} with $\kappa=0$ on $\mathbb{T}^3\times[0,\infty)$ with radius of convergence at least $\tau=\tau(t)$ for all $t\in[0,\infty)$, where $\tau$ is a decreasing function satisfying
\begin{align*}
\tau(t)\ge\tau_0e^{-C\left(\|e^{\tau_0\Lambda^\frac{1}{s}}\theta_0\|_{L^2}+2\|e^{\tau_0\Lambda^\frac{1}{s}}S\|_{L^2}\right)t}.
\end{align*}
Here $C>0$ is a constant which depends on $\nu$ but independent of $t$. For the singular case when $\nu=0$, there exists $\bar{\tau}\in(0,\tau_0]$, $\bar{T}>0$ and a unique Gevrey-class $s$ solution $\theta^0$ to \eqref{MG equation}-\eqref{MG Fourier symbol} for $\kappa=0$ defined on $\mathbb{T}^3\times[0,\bar{T}]$ with radius of convergence at least $\bar{\tau}$.
\end{thm}

\begin{thm}[Convergence of solutions as $\nu\rightarrow0$ for the MG$^\nu$ equations]\label{Convergence of solutions as nu goes to 0 non-viscid MG thm}
Fix $s\ge1$, $r>3$ and $K_0>0$. Let $\theta_0$ and $S$ be of Gevrey-class $s$ with radius of convergence $\tau_0>0$ and satisfy the assumptions given in Theorem~\ref{Gevrey-class global well-posedness MG non-viscid thm}. If $\theta^\nu$ and $\theta^0$ are Gevrey-class $s$ solutions to \eqref{MG equation}-\eqref{MG Fourier symbol} with $\kappa=0$ for $\nu>0$ and $\nu=0$ respectively with initial datum $\theta_0$ on $\mathbb{T}^3\times[0,\bar{T}]$ with radius of convergence at least $\bar{\tau}$ as described in Theorem~\ref{Gevrey-class global well-posedness MG non-viscid thm}, then there exists $T<\bar{T}$ and $\tau=\tau(t)<\bar{\tau}$ such that, for $t\in[0,T]$, we have
\begin{align*}
\lim_{\nu\rightarrow0}\|(\Lambda^re^{\tau\Lambda^\frac{1}{s}}\theta^{\nu}-\Lambda^re^{\tau\Lambda^\frac{1}{s}}\theta^0)(\cdot,t)\|_{L^2}=0.
\end{align*}
\end{thm}

\begin{thm}[Well-posedness in Sobolev space for the IPMB$^\nu$ equations]\label{Well-posedness in Sobolev space IPMB thm}
Let $\theta_0\in W^{s,2}$ for $s\ge0$. Then for each $\nu>0$, we have:
\begin{itemize}
\item if $s=0$, there exists unique global weak solution to \eqref{IPMB constitutive law}-\eqref{Fourier multiplier symbol IPMB} such that
\begin{align*}
\theta^\nu&\in BC((0,\infty);L^2),\\
u^\nu&\in C((0,\infty);W^{2,2}).
\end{align*}
In particular, $\theta^\nu(\cdot,t)\rightarrow\theta_0$ weakly in $L^2$ as $t\rightarrow0^+$.
\item if $s>0$, there exists a unique global-in-time solution $\theta^\nu$ to \eqref{IPMB constitutive law}-\eqref{Fourier multiplier symbol IPMB} such that $\theta^\nu(\cdot,t)\in W^{s,2}$ for all $t\ge0$. Furthermore, for $s=1$, we have the following single exponential growth in time on $\|\nabla\theta^\nu(\cdot,t)\|_{L^2}$:
\begin{align*}
\|\nabla\theta^\nu(\cdot,t)\|_{L^2}\le C \|\nabla\theta_0\|_{L^2}\exp\left(Ct\|\theta_0\|_{W^{1,2}}\right),
\end{align*}
where $C>0$ is a constant which depends only on some dimensional constants.
\end{itemize}
\end{thm}

\begin{thm}[Gevrey-class global well-posedness for the IPMB$^\nu$ equations]\label{Gevrey-class global well-posedness IPMB thm}
Fix $s\ge1$. Let $\theta_0$ be of Gevrey-class $s$ with radius of convergence $\tau_0>0$. Then for each $\nu>0$, there exists a unique Gevrey-class $s$ solution $\theta^\nu$ to \eqref{IPMB constitutive law}-\eqref{Fourier multiplier symbol IPMB} on $\mathbb{T}^2\times[0,\infty)$ with radius of convergence at least $\tau=\tau(t)$ for all $t\in[0,\infty)$, where $\tau$ is a decreasing function satisfying
\begin{align*}
\tau(t)\ge\tau_0e^{-Ct\|e^{\tau_0\Lambda^\frac{1}{s}}\theta_0\|_{L^2}}.
\end{align*}
Here $C>0$ is a constant which depends on $\nu$ but independent of $t$.
\end{thm}

\begin{thm}[Local well-posedness and convergence of solutions in Sobolev space for the IPMB$^\nu$ equations]\label{local-in-time existence and convergence theorem for IPMB thm}
Fix $s>2$ and assume that $\theta_0\in H^s(\mathbb{T}^2)$ has zero-mean on $\mathbb{T}^2$. Then there exists a positive time $T$ and a unique smooth solution $\theta^0$ to \eqref{IPMB constitutive law}-\eqref{Fourier multiplier symbol IPMB} with $\nu=0$ such that
$$\theta^0\in L^\infty(0,T;H^s(\mathbb{T}^2)).$$
Moreover, for $t\in[0,T]$, we have
\begin{align*}
\lim_{\nu\rightarrow0}\|(\theta^{\nu}-\theta^0)(\cdot,t)\|_{H^{s-1}}=0.
\end{align*}
\end{thm}

\section{Fractionally diffusive active scalar equations}\label{s:Fractionally diffusive active scalar equations}

In this section, we investigate the properties of the family of active scalar equations \eqref{abstract active scalar eqn}-\eqref{zero mean assumption} in the context of the fractional Laplacian. The results can be applied to the modified surface quasi-geostrophic (SQG$^{\kappa,\nu}$) equation \eqref{surface quasi-geostrophic general SQG}-\eqref{Fourier multiplier symbol SQG} and MG$^\nu$ equation \eqref{MG equation}-\eqref{MG Fourier symbol}.

\subsection{Existence and convergence of $H^s$-solutions when $\nu >0$}\label{s:Existence and convergence of Hs-solutions when nu>0}

When the parameter $\nu$ is taken to be positive, the ensuing smoothing properties of $T_{ij}^{\nu}$ permits existence and convergence in Sobolev space $H^s$ as $\kappa$ goes to zero. The global-in-time existence theorem is given as follows:

\begin{thm}[Friedlander and Suen \cite{FS20}]\label{global-in-time well-posedness in Sobolev fractional thm}
Fix $\nu>0$, $s\ge0$ and $\gamma\in(0,2]$, and let $\theta_0\in H^s$ and $S\in H^s\cap L^\infty$ be given. 
\begin{itemize}
\item For any $\kappa>0$, there exists a global-in-time solution to \eqref{abstract active scalar eqn} such that 
\begin{align*}
\theta^\kappa\in C([0,\infty);H^s)\cap L^2([0,\infty);H^{s+\frac{\gamma}{2}}).
\end{align*}
\item For $\kappa=0$, if we further assume that $\theta_0\in L^\infty$, then exists a global-in-time solution to \eqref{abstract active scalar eqn}-\eqref{zero mean assumption}  such that $\theta^0(\cdot,t)\in H^s$ for all $t\ge0$.
\end{itemize}
\end{thm}

In view of the case when $\kappa>0$, the most subtle part for proving Theorem~\ref{global-in-time well-posedness in Sobolev fractional thm} is to estimate the $L^\infty$-norm of $\theta^\kappa(\cdot,t)$ when $\theta_0$ is {\it not} necessarily in $L^\infty$. In achieving our goal, we apply De Giorgi iteration method (see \cite[Lemma~4.5]{FS20}) and obtain
\begin{align*}
\|\theta^\kappa(t)\|_{L^\infty}\le C\Big[\Big(\frac{1}{\kappa t}\Big)^{\frac{d+1-\gamma}{2\gamma}}\Big(\|\theta_0\|_{L^2}+\frac{\|S\|_{L^2}}{c_0^\frac{1}{2}\kappa^\frac{1}{2}}\Big)+\|S\|_{L^\infty}\Big],
\end{align*}
for some constant $C=C(d)>0$ which only depends on the dimension $d$. Once Theorem~\ref{global-in-time well-posedness in Sobolev fractional thm} is proved, we can show the convergence of $H^s$ solutions which are summarised as follows:
\begin{thm}[Friedlander and Suen \cite{FS20}]\label{Hs convergence fractional thm}
Let $\nu>0$ and $\gamma\in(0,2]$ be given in \eqref{abstract active scalar eqn}, and let $\theta_0,S\in C^\infty$ be the initial datum and forcing term respectively which satisfy \eqref{zero mean assumption on data and forcing}. If $\theta^\kappa$ and $\theta^0$ are smooth solutions to \eqref{abstract active scalar eqn}-\eqref{zero mean assumption}  for $\kappa>0$ and $\kappa=0$ respectively, then 
\begin{align*}
\lim_{\kappa\rightarrow0}\|(\theta^\kappa-\theta^0)(\cdot,t)\|_{H^s}=0,
\end{align*}
for all $s\ge0$ and $t\ge0$.
\end{thm}

\subsection{Existence and convergence of Gevrey-class solutions when $\nu=0$}\label{s:Existence and convergence of Gevrey-class solutions}

In contrast to the case for $\nu>0$, when the parameter $\nu$ is set to zero, condition (A2) implies that $\partial_{x_i} T_{ij}^{\nu}$ is a singular operator. In this case the existence and convergence results for \eqref{abstract active scalar eqn}-\eqref{zero mean assumption}  are restricted to Gevrey-class solutions which are summarised in the following theorem:

\begin{thm}[Friedlander and Suen \cite{FS20}]\label{Local-in-time existence of Gevrey-class solutions fractional thm}
Let $\kappa\ge0$ and $\gamma\in(0,2]$ be fixed, and let $\theta_0$ and $S$ be the initial datum and forcing term respectively. Fix $s\ge1$ and $K_0>0$. Suppose $\theta_0$ and $S$ both belong to Gevrey-class $s$ with radius of convergence $\tau_0>0$ and
\begin{align*}
\|\Lambda^re^{\tau_0\Lambda^\frac{1}{s}}\theta^\kappa(\cdot,0)\|_{L^2}\le K_0,\qquad\|\Lambda^re^{\tau_0\Lambda^\frac{1}{s}}S\|_{L^2}\le K_0,
\end{align*}
where $r>\max\{\frac{d}{2}+\frac{3}{2},\frac{d}{2}+1+\frac{3}{2s}\}$. Then there exists $T_*=T_*(\tau_0, K_0) > 0$ and a unique Gevrey-class $s$ solution on $[0, T_*)$ to the initial value problem associated to \eqref{abstract active scalar eqn}-\eqref{zero mean assumption}  with $\nu=0$. Furthermore, if $\theta^{\kappa}$, $\theta^{0}$ are Gevrey-class $s$ solutions to \eqref{abstract active scalar eqn}-\eqref{zero mean assumption}  with $\nu=0$ for $\kappa>0$ and $\kappa=0$ respectively with initial datum $\theta_0$ on $\mathbb{T}^d\times[0,\bar{T}]$ with radius of convergence at least $\bar{\tau}$, then there exists $T\le\bar{T}$ and $\tau=\tau(t)<\bar{\tau}$ such that, for $t\in[0,T]$, we have:
\begin{align*}
\lim_{\kappa\to0}\|(\Lambda^re^{\tau\Lambda^\frac{1}{s}}\theta^{\kappa}-\Lambda^re^{\tau\Lambda^\frac{1}{s}}\theta^{0})(\cdot,t)\|_{L^2}=0.
\end{align*}
\end{thm}

For $\kappa>0$ and $\gamma\in[1,2]$, under a smallness assumption on the initial data, it can be proved that the Gevrey-class solutions obtained in Theorem~\ref{Local-in-time existence of Gevrey-class solutions fractional thm} exist for all time:

 \begin{thm}[Friedlander and Suen \cite{FS20}]\label{Global-in-time existence of Gevrey-class solutions fractional thm}
Let $\kappa>0$ and $\gamma\in[1,2]$, and suppose that both $\theta_0$ and $S$ belong to some Gevrey-class $s$ with $s\ge1$. There exists $\varepsilon>0$ depending on $\kappa$ such that, if $\theta_0$ and $S$ satisfy
\begin{equation}\label{smallness assumption on theta0 and S 1}
\|\theta_0\|_{L^2}^\beta\|\theta_0\|_{H^\alpha}^{1-\beta}+\|\theta_0\|_{L^2}^\beta\|S\|^{1-\beta}_{L^\infty(0,\infty;H^{\alpha-\gamma})}\le\varepsilon,
\end{equation}
and
\begin{equation*}
\|\Lambda^\alpha\theta_0\|^2_{L^2}+\frac{2}{\kappa^2}\|S\|^2_{H^{\alpha-\gamma}}\le\varepsilon,
\end{equation*} 
where $\alpha>\frac{1}{2}(d+2)+(1-\gamma)$ and $\beta=1-\frac{1}{\alpha}\Big[\frac{1}{2}(d+2)+(1-\gamma)\Big]$, then the local-in-time Gevrey-class $s$ solution $\theta^\kappa$ as claimed by Theorem~\ref{Local-in-time existence of Gevrey-class solutions fractional thm} can be extended to all time.
\end{thm}

As a by-product of Theorem~\ref{Global-in-time existence of Gevrey-class solutions fractional thm}, for the case when $S\in H^{s-\gamma}(\mathbb{T}^d)$ and $\theta_0\in H^\alpha(\mathbb{T}^d)$ with $\gamma\in[1,2]$ and $\alpha>\frac{1}{2}(d+2)+(1-\gamma)$, under the smallness assumption \eqref{smallness assumption on theta0 and S 1}, the system \eqref{abstract active scalar eqn}-\eqref{zero mean assumption} possesses a global-in-time $H^\alpha$ solution:

\begin{thm}[Friedlander and Suen \cite{FS20}]\label{Global-in-time existence of Sobolev solutions thm}
Let $\kappa>0$, $\gamma\in[1,2]$ and $S\in H^{s-\gamma}(\mathbb{T}^d)$, and let $\theta_0\in H^\alpha(\mathbb{T}^d)$ have zero mean on $\mathbb{T}^d$, where $\alpha>\frac{d+2}{2} + (1-\gamma)$.
There exists a small enough constant $\varepsilon>0$ depending on $\kappa$, such that if $\theta_0$ satisfies \eqref{smallness assumption on theta0 and S 1}, then there exists a unique global-in-time $H^\alpha$-solution to \eqref{abstract active scalar eqn}-\eqref{zero mean assumption} with $\nu=0$. In particular, for all $t>0$, we have the following bound on $\theta^\kappa$:
\begin{equation*}
\|\Lambda^\alpha\theta^\kappa(\cdot,t)\|^2_{L^2}\le \|\Lambda^\alpha\theta_0\|^2_{L^2}+\frac{2}{\kappa^2}\|S\|^2_{H^{\alpha-\gamma}}.
\end{equation*}
\end{thm}

The proofs of Theorem~\ref{Local-in-time existence of Gevrey-class solutions fractional thm}-\ref{Global-in-time existence of Sobolev solutions thm} can be found in \cite{FS20}. We point out that all the abstract results obtained in Theorem~\ref{Local-in-time existence of Gevrey-class solutions fractional thm}-\ref{Global-in-time existence of Sobolev solutions thm} can be applied to the critical SQG equation, which is a special example of \eqref{abstract active scalar eqn}-\eqref{zero mean assumption} with $\nu = 0$ and $\gamma = 1$.

\subsection{Long time behaviour for solutions when $\nu > 0$ and $\kappa > 0$}\label{s:Long time behaviour for solutions when nu > 0 and kappa > 0}

In this subsection, we study the long time behaviour for solutions to the active scalar equations \eqref{abstract active scalar eqn} when $\nu>0$ and $\kappa>0$. Based on the global-in-time existence results established in Theorem~\ref{global-in-time well-posedness in Sobolev fractional thm}, for fixed $\nu>0$ and $\kappa>0$, we can define a solution operator $\pin(t)$ for the initial value problem \eqref{abstract active scalar eqn} via
\begin{align}\label{def of solution map nu>0 and kappa>0}
\pin(t): H^1\to H^1,\qquad \pin(t)\theta_0=\theta(\cdot,t),\qquad t\ge0.
\end{align}
We study the long-time dynamics of $\pin(t)$ on the phase space $H^1$. Specifically, we establish the existence of global attractors for $\pin(t)$ in $H^1$ and address some properties for the attractors. The following theorem first gives the existence of global attractors:

\begin{thm}[Friedlander and Suen \cite{FS20}]\label{Existence of H1 global attractor thm}
Let $S\in L^\infty\cap H^1$. For $\nu$, $\kappa>0$ and $\gamma\in(0,2]$, the solution map $\pin(t):H^1\to H^1$ associated to \eqref{abstract active scalar eqn} possesses a unique global attractor $\Gg$. Moreover, there exists $M_{\Gg}$ which depends only on $\nu$, $\kappa$, $\gamma$, $\|S\|_{L^\infty\cap H^1}$ and universal constants, such that if $\theta_0\in\Gg$, we have that
\begin{align}\label{uniform bound on theta from the attractor} 
\|\theta(\cdot,t)\|_{H^{1+\frac{\gamma}{2}}}\le M_{\Gg},\qquad\forall t\ge0,
\end{align}
and
\begin{align}\label{uniform bound on time integral of theta from the attractor}
\frac{1}{T}\int_t^{t+T}\|\theta(\cdot,\tau)\|_{H^{1+\gamma}}d\tau\le M_{\Gg},\qquad \mbox{$\forall t\ge0$ and $T>0$,}
\end{align}
where $\theta(\cdot,t)=\pin(t)\theta_0$.
\end{thm}

Details of the proof of Theorem~\ref{Existence of H1 global attractor thm} can be found in \cite[Subsection~6.1]{FS20}, we also refer to \cite{CZV16} for the case of SQG equations. The steps of proof can be outlined as follows:
\begin{itemize}
\item[I.] By standard energy method (see for example \cite{CTV14a} for the case $\gamma=1$), one can show that $\theta^\kappa$ satisfies
\begin{align}\label{energy inequality kappa>0}
\|\theta^\kappa(\cdot,t)\|^2_{L^2}+\kappa\int_0^t\|\Lambda^\frac{\gamma}{2}\theta^\kappa(\cdot,\tau)\|^2_{L^2}d\tau\le\|\theta_0\|^2_{L^2}+\frac{t}{c_0\kappa}\|S\|^2_{L^2}, \qquad\forall t\ge0,
\end{align}
where $c_0>0$ is a universal constant which depends only on the dimension $d$. Moreover, by \cite[Lemma~6.2]{FS20}, the set 
\begin{align*}
B_{\infty}=\left\{\phi\in L^\infty\cap H^1:\|\phi\|_{L^\infty}\le\frac{2}{c_0\kappa}\|S\|_{L^\infty}\right\}
\end{align*}
is an absorbing set for $\pin(t)$ and
\begin{align}\label{bound on solution map in L infty}
\sup_{t\ge0}\sup_{\theta_0\in B_{\infty}}\|\pin(t)\theta_0\|_{L^\infty}\le\frac{3}{c_0\kappa}\|S\|_{L^\infty}.
\end{align}
\item[II.] Next by \cite[Lemma~6.3]{FS20},  we obtain the necessary {\it a priori} estimate in $C^{\alpha}$-space with some appropriate exponent $\alpha\in(0,1)$. Furthermore, as pointed out in \cite{CTV14a}, we can see that the solutions to \eqref{abstract active scalar eqn} emerging from data in a bounded subset of $H^1$ are absorbed in finite time by $B_{\infty}$. Hence if $\theta_0\in H^1\cap L^\infty$ and fix $\nu$, $\kappa>0$, then there exists $\alpha\in(0,\frac{\gamma}{2\gamma+2}]$ which depends on $\|\theta_0\|_{L^\infty}$, $\|S\|_{L^\infty}$, $\nu$, $\kappa$, $\gamma$ such that
\begin{align}\label{C alpha bounds}
\|\theta(\cdot,t)\|_{C\alpha}\le C(\K+\bK),\qquad\forall t\ge t_{\alpha}:=\frac{2\gamma(1-\alpha)}{2+\gamma},
\end{align}
where $C>0$ is a positive constant, $\K$ and $\bK$ are given respectively by
\begin{align*}
\K:=\|\theta_0\|_{L^\infty}+\frac{1}{c_0\kappa}\|S\|_{L^\infty},\qquad \bK:=\kappa^{-1}\K^{\gamma}+\kappa^{-\frac{1}{4}}\|S\|^{\frac{2+\gamma}{2(1+\gamma)}}_{L^\infty}\K^\frac{\gamma}{4}+\kappa^{-1}\K^{\gamma+\frac{2+\gamma}{2\gamma}}.
\end{align*}
With the help of \eqref{C alpha bounds}, we obtain the following result which can be regarded as an improvement of the regularity of the absorbing set $B_{\infty}$, namely there exists $\alpha\in(0,\frac{\gamma}{2\gamma+2}]$ and a constant $C_{\alpha}=C_{\alpha}(\|S\|_{L^\infty},\alpha,\nu,\kappa,\gamma,\K,\bK)\ge1$ such that the set
\begin{align*}
B_{\alpha}=\left\{\phi\in C^{\alpha}\cap H^1:\|\phi\|_{C^\alpha}\le C_{\alpha}\right\}
\end{align*}
is an absorbing set for $\pin(t)$. 
\item[III.] As in \cite[Lemma~6.7]{FS20}, by establishing an {\it a priori} estimate for initial data in $H^1\cap C^\alpha$, we can show that there exists a bounded absorbing set for $\pin(t)$ in $H^1$. More precisely, there exists $\alpha\in(0,\frac{\gamma}{2\gamma+2}]$ and a constant $R_1=R_1(\|S\|_{L^\infty\cap H^1},\alpha,\nu,\kappa,\gamma)\ge1$ such that the set
\begin{align*}
B_1=\{\phi\in C^{\alpha}\cap H^1:\|\phi\|^2_{H^1}+\|\phi\|^2_{C^\alpha}\le R^2_{1}\}
\end{align*}
is an absorbing set for $\pin(t)$. Moreover, we have
\begin{align}\label{bound on time integral on theta in H 1+gamma/2}
\sup_{t\ge0}\sup_{\theta_0\in B_1}\left[\|\pin(t)\theta_0\|^2_{H^1}+\|\pin(t)\theta_0\|^2_{C^\alpha}+\int_t^{t+1}\|\pin(\tau)\theta_0\|^2_{H^{1+\frac{\gamma}{2}}}d\tau\right]\le 2R_1^2.
\end{align}
The bound \eqref{bound on time integral on theta in H 1+gamma/2} turns out to be crucial for improving the regularity of the absorbing set $B_1$ to $H^{1+\frac{\gamma}{2}}$, which allows us to obtain an absorbing set $B_{1+\frac{\gamma}{2}}$ for $\pin(t)$ given by
\begin{align*}
B_{1+\frac{\gamma}{2}}=\left\{\phi\in H^{1+\frac{\gamma}{2}}:\|\phi\|_{H^{1+\frac{\gamma}{2}}}\le R_{1+\frac{\gamma}{2}}\right\}
\end{align*}
for some constant $R_{1+\frac{\gamma}{2}}\ge1$ which depends on $\|S\|_{L\infty\cap H^1}$, $\nu$, $\kappa$, $\gamma$. The existence and regularity of the global attractor claimed by Theorem~\ref{Existence of H1 global attractor thm} now follows by applying the argument given in \cite[Proposition 8]{CCP12} and the bound \eqref{bound on time integral on theta in H 1+gamma/2}.
\end{itemize}

After we have obtained the global attractors as described in Theorem~\ref{Existence of H1 global attractor thm}, we prove some additional properties on the attractors under the assumption that $\gamma\in[1,2]$ (also refer to \cite[Subsection~6.2]{FS20} for details):

\begin{thm}[Friedlander and Suen \cite{FS20}]\label{Further properties of H1 global attractor thm}
Let $S\in L^\infty\cap H^1$. For $\nu$, $\kappa>0$, assume that the exponent $\gamma\in[1,2]$. Then the global attractor $\Gg$ of $\pin(t)$ further enjoys the following properties:
\begin{itemize}
\item $\Gg$ is fully invariant, namely
\begin{align*}
\pin(t)\Gg=\Gg,\qquad \forall t\ge0.
\end{align*}
\item $\Gg$ is maximal in the class of $H^1$-bounded invariant sets.
\item $\Gg$ has finite fractal dimension.
\end{itemize}
\end{thm}

To prove the invariance and the maximality of the attractor $\Gg$, we observe that the solution map $\pin(t)$ is indeed continuous in the $H^1$-topology, in other words for every $t>0$, the solution map $\pin(t):B_{1+\frac{\gamma}{2}}\to H^1$ is continuous in the topology of $H^1$. The key ingredient for the proof of continuity is the bound
\begin{align*}
\|u\|_{L^\infty}\le C_\nu\|\Lambda\theta\|_{L^2},
\end{align*}
where $C_\nu>0$ is a constant depending on $\nu>0$, and such bound comes from the condition (A3) and the assumption that $d=2$ or 3. Following the argument given in \cite[Proposition 5.5]{CTV14a} and using the log-convexity method introduced by \cite{AN67}, we can also prove that the solution map $\pin$ is injective on the absorbing set $B_{1+\frac{\gamma}{2}}$. Hence by applying \cite[Proposition~6.4]{CZV16}, we can obtain the invariance and the maximality of the attractor $\Gg$ stated in Theorem~\ref{Further properties of H1 global attractor thm}. 

It remains to address the fractal dimensions for the global attractors $\Gg$. In order to prove that $\dimf (\Gg)$ is finite, we need to show that the solution map $\pin$ is {\it uniform differentiable} (refer to \cite[Definition~6.18]{FS20} for the definition for being uniform differentiable). And by \cite[Lemma~6.21]{FS20}, the large-dimensional volume elements which are carried by the flow of $\pin(t)\theta_0$, with $\theta_0\in\Gg$, actually have exponential decay in time. The argument in \cite[pp. 115--130, and Chapter 14]{CF88} can then be applied which shows that $\dimf(\Gg)$ is finite.

\subsection{Applications to magneto-geostrophic equations}

We now apply our results claimed in subsections~\ref{s:Existence and convergence of Hs-solutions when nu>0}-\ref{s:Long time behaviour for solutions when nu > 0 and kappa > 0} to MG$^\nu$ equation \eqref{MG equation}-\eqref{MG Fourier symbol}, which can be summarised in the following theorems (see also \cite[Subsection~7.1]{FS20}. We point out that, Theorem~\ref{Hs convergence MG} strengthens and generalises the results obtained in \cite{FS15} in which the authors showed weak convergence as $\kappa\to0$.

\begin{thm}[$H^s$-convergence as $\kappa\rightarrow0$ for MG$^\nu$ equations]\label{Hs convergence MG}
Let $\nu>0$ be given as in \eqref{MG equation}, and let $\theta_0,S\in C^\infty$ be the initial datum and forcing term respectively. If $\theta^\kappa$ and $\theta^0$ are smooth solutions to \eqref{MG equation}-\eqref{MG Fourier symbol} for $\kappa>0$ and $\kappa=0$ respectively, then 
\begin{align*} 
\lim_{\kappa\rightarrow0}\|(\theta^\kappa-\theta^0)(\cdot,t)\|_{H^s}=0,
\end{align*}
for all $s\ge0$ and $t\ge0$.
\end{thm}

\begin{thm}[Gevrey-class convergence as $\kappa\rightarrow0$ for MG equations]\label{Analytic convergence kappa MG}
Let $\nu=0$ be given as in \eqref{MG equation}, and let $\theta_0,S$ the initial datum and forcing term respectively. Suppose that both $\theta_0$ and $S$ belong to some Gevrey-class $s$ with $s\ge1$. Then if $\theta^{\kappa}$, $\theta^{0}$ are Gevrey-class $s$ solutions to \eqref{MG equation}-\eqref{MG Fourier symbol} for $\kappa>0$ and $\kappa=0$ respectively with initial datum $\theta_0$ and with radius of convergence at least $\bar{\tau}$, then there exists $T\le\bar{T}$ and $\tau=\tau(t)<\bar{\tau}$ such that, for $t\in[0,T]$, we have:
\begin{align*}
\lim_{\kappa\to0}\|(\Lambda^re^{\tau\Lambda^\frac{1}{s}}\theta^{\kappa}-\Lambda^re^{\tau\Lambda^\frac{1}{s}}\theta^{0})(\cdot,t)\|_{L^2}=0.
\end{align*}
\end{thm}

\begin{thm}[Existence of global attractors for MG$^\nu$ equations]\label{existence of global attractor MG theorem}
Let $S\in L^\infty\cap H^1$. For $\nu$, $\kappa>0$, let $\pin(t)$ be solution operator for the initial value problem \eqref{MG equation} via \eqref{def of solution map nu>0 and kappa>0}. Then the solution map $\pin(t):H^1\to H^1$ associated to \eqref{MG equation}-\eqref{MG Fourier symbol} possesses a unique global attractor $\Gg$ for all $\nu>0$. In particular, for each $\nu>0$, the global attractor $\Gg$ of $\pin(t)$ enjoys the following properties:
\begin{itemize}
\item $\Gg$ is fully invariant, namely
\begin{align*}
\pin(t)\Gg=\Gg,\qquad \forall t\ge0.
\end{align*}
\item $\Gg$ is maximal in the class of $H^1$-bounded invariant sets.
\item $\Gg$ has finite fractal dimension.
\end{itemize}
\end{thm}

We recall from subsection~\ref{s:The MG equations and existence of compact global attractor diffusive} that there exists a compact global attractor $\A$ in $L^2(\mathbb{T}^3)$ for the MG$^0$ equations, namely the equations \eqref{MG equation}-\eqref{MG Fourier symbol} when $\kappa>0$, $\nu=0$ and $S\in L^\infty\cap H^1$. When $\nu$ is varying, we relate the global attractors $\Gg$ with $\A$ and further obtain the following theorem:

\begin{thm}[Friedlander and Suen \cite{FS20}]\label{varying nu theorem}
Let $\kappa>0$ be fixed in \eqref{MG equation}. Then we have:
\begin{enumerate}
\item If $\Gg$ are the global attractors for the MG$^\nu$ equations \eqref{MG equation}-\eqref{MG Fourier symbol} as obtained by Theorem~\ref{existence of global attractor MG theorem}, then $\Gg$ and $\A$ satisfy
\begin{align}\label{upper semi-continuity at nu=0}
\mbox{$\dis\sup_{\phi\in\Gg}\inf_{\psi\in\A}\|\phi-\psi\|_{L^2}\rightarrow0$ as $\nu\rightarrow0$.}
\end{align}
\item Let $\nu^*>\nu_*>0$ be arbitrary. For each $\nu_0\in[\nu_*,\nu^*]$, the collection $\dis\{\Gg\}_{\nu\in[\nu_*,\nu^*]}$ is {\it upper semicontinuous} at $\nu_0$ in the following sense: 
\begin{align}\label{upper semi-continuity at fixed nu>0}
\mbox{$\dis\sup_{\phi\in\Gg}\inf_{\psi\in\Ggo}\|\phi-\psi\|_{H^1}\rightarrow0$ as $\nu\rightarrow\nu_0$.}
\end{align}
\end{enumerate}
\end{thm}

To prove the convergence \eqref{upper semi-continuity at nu=0}, we recall from Theorem~\ref{existence of MG attractor thm} that for $\kappa>0$, $\nu\in[0,1]$ and $S\in L^\infty\cap H^2$, there exists global attractor $\An$ in $L^2$ generated by the solution map $\pino$ via
\begin{align*}
\pino(t): L^2\to L^2,\qquad \pino(t)\theta_0=\theta(\cdot,t),\qquad t\ge0,
\end{align*}
and $\An$ is upper semicontinuous at $\nu=0$. Since $\pino\Big|_{H^1}=\pin$ and $\Gg\subset \An$ for all $\nu\in(0,1]$, the convergence \eqref{upper semi-continuity at nu=0} follows immediately from \eqref{uc_0}. On the other hand, to prove the convergence \eqref{upper semi-continuity at fixed nu>0}, we write $\I=[\nu_*,\nu^*]$ and show that
\begin{itemize}
\item[I.] there is a compact subset $\U$ of $H^1$ such that $\Gg\subset\U$ for every $\nu\in\I$; and
\item[II.] for $t > 0$, $\pin\theta_0$ is continuous in $\I$, uniformly for $\theta_0$ in compact subsets of $H^1$.
\end{itemize}
The key for showing Step I. and Step II. is the following bound, namely for any $\nu\in\I$, $s\in[0,2]$ and $f\in L^p$ with $p>1$, we have
\begin{align}\label{two order smoothing for u when nu inside I*}
\|\Lambda^s u^{\nu}[f]\|_{L^p}\le C_*\|f\|_{L^p},
\end{align}
where $C_*$ is a positive constant which depends only on $p$, $\nu_*$ and $\nu^*$. The bound can be used for obtaining a bounded set $B_2$ in $H^2$ given by 
\begin{align*}
B_{2}=\left\{\phi\in H^{2}:\|\phi\|_{H^{2}}\le R_{2}\right\}
\end{align*}
where $R_{2}\ge1$ is a constant which depends only on $\nu_*$, $\nu^*$, $\kappa$, $\|S\|_{L^\infty\cap H^1}$, and $B_2$ enjoys the following properties:
\begin{itemize}
\item $B_{2}$ is a compact set in $H^1$ which depends only on $\nu_*$, $\nu^*$, $\kappa$, $\|S\|_{L^\infty\cap H^1}$;
\item $\Gg\subset B_{2}$ for all $\nu\in\I$.
\end{itemize}
For the uniform continuity stated in Step II., with the help of the bound \eqref{two order smoothing for u when nu inside I*}, we can obtain the necessary $H^1$-estimates \cite[Lemma~7.11]{FS20}:
Define $\U=\{\phi\in H^1:\|\phi\|^2_{H^1}\le R_{\U}\}$ where $R_{\U}>0$, then for any $\theta_0\in\U$ and $\nu\in\I$, if $\theta^\nu(t)=\pin(t)\theta_0$, then $\theta^\nu(t)$ satisfies
\begin{align}\label{H1 estimate in compact set of H1}
\sup_{0\le \tau\le t}\|\theta^\nu(\cdot,\tau)\|^2_{H^1}+\int_0^t\|\theta^\nu(\cdot,\tau)\|^2_{H^2}d\tau\le M_*(t),\qquad\forall t>0,
\end{align}
where $M_*(t)$ is a positive function in $t$ which depends only on $t$, $\kappa$, $\nu_*$, $\nu^*$, $\|S\|_{H^1}$ and $R_{\U}$. The uniform continuity stated in Step II. then follows by energy-type estimates on the difference $\theta^{\nu_1}-\theta^{\nu_2}$ with $\nu_1$, $\nu_2\in\I$, which completes the proof of Theorem~\ref{varying nu theorem}.


\section*{Acknowledgment} S. Friedlander is supported by NSF DMS-1613135 and A. Suen is supported by Hong Kong General Research Fund (GRF) grant project number 18300720.

\end{document}